# Fourier series (based) multiscale method for computational analysis in science and engineering:

# V. Fourier series multiscale solution for elastic bending of Reissner plates on Pasternak foundations

Weiming Sun[a*+] and Zimao Zhang[b]

**Abstract:** Fourier series multiscale method, a concise and efficient analytical approach for multiscale computation, will be developed out of this series of papers. In the fifth paper, the usual structural analysis of plates on an elastic foundation is extended to a thorough multiscale analysis for a system of a fourth order linear differential equation (for transverse displacement of the plate) and a second order linear differential equation (for the stress function), where general boundary conditions and a wide spectrum of model parameters are prescribed. For this purpose, the solution function each is expressed as a linear combination of the corner function, the two boundary functions and the internal function, to ensure the series expression obtained uniformly convergent and termwise differentiable up to fourth (or second) order. Meanwhile, the sum of the corner function and the internal function corresponds to the particular solution, and the two boundary functions correspond to the general solutions which satisfy the homogeneous form of the equation. Since the general solutions have appropriately interpreted the meaning of the differential equation, the spatial characteristics of the solution of the equation are expected to be better captured in separate directions. With the corner function, the two boundary functions and the internal function selected specifically as polynomials, one-dimensional half-range Fourier series along the $y$ (or $x$)-direction, and two-dimensional half-range Fourier series, the Fourier series multiscale solution of the bending problem of a Reissner plate on the Pasternak foundation is derived. And then the convergence characteristics of the Fourier series multiscale solution are investigated with numerical examples, and the multiscale characteristics of the bending problem of a Reissner plate on the Pasternak foundation are demonstrated for a wide spectrum of model parameters.

---

[a]Department of Mathematics and Big Data, School of Artificial Intelligence, Jianghan University, Wuhan, 430056, China

[b]Department of Mechanics, School of Civil Engineering, Beijing Jiaotong University, Beijing, 100044, China

*Correspondence to: Weiming Sun, Department of Mathematics and Big Data, School of Artificial Intelligence, Jianghan University, Wuhan, 430056, China

[+] E-mail: xuxinenglish@hust.edu.cn

# 1. Introduction

Nowadays, heavy vehicles and large aircrafts have started to become common means of transport and meanwhile the aerospace technology have been developed rapidly. As a result, the thicknesses of reinforced concrete pavements of various types of engineering structures, such as airport runways, highways, and space launch sites increase continuously. For example, the thicknesses of concrete pavements of heavy duty runways and parking yards of container vehicles are all greater than 50 cm. However, during the design and construction process of these concrete pavements, traditional thin plate bending theory is still employed for the structural analysis of the foundation-plate systems. And then, the thicknesses of concrete pavements are artificially increased to ensure service lives of the structures, which may result in tremendous wastes of resources and money. On the other hand, if the concrete pavements are too thin, the service lives of the structures will be shortened greatly. and worse, serious road accidents arise. In these circumstances, important launch missions are delayed or traffic jams are brought about. Therefore, thick plate theory is desired for more accurate analysis of concrete pavements of highways, airfields and space launch sites, and of the floor system of industrial yards.

Unfortunately, the complexity of the analytical formulations of thick plates on elastic foundation limited the number of available analytical solutions. As the result, several numerical methods have been used by researchers to solve the plate bending problem. These numerical methods include finite difference method [1], finite element method [2, 3], boundary element method [4, 5], meshless method [6], numerical manifold method [7], finite element method of lines [8], method of fundamental solution (MFS) [9], etc.

As a general analytical method for differential equations, the Fourier series method not only is of great theoretical significance, but also is of important practical value for solving the bending problem of thick plates on elastic foundations. For example, based on the linear distribution of bending normal stress through the thickness of rectangular plate, Henwood described the linear elastic behavior of a thick plate on the Winkler foundation as a system of eight first order partial differential equations, solved this problem with the general analytical Fourier series method, and obtained the corresponding series solution [10]. For thick plates with four free edges, he made a comparative analysis of the convergence characteristics, computational accuracy and computational efficiency of the series solution under various types of loading conditions. It is shown that the Fourier series solution is of good convergence and high computational accuracy, and can satisfy the boundary conditions (except that near corners), and is of higher computational efficiency than that of the finite difference method and the finite element method. Shi and Yao investigated the problem of rectangular thick plates with four free edges on an elastic foundation by the Fourier series based superposition method, where the Reissner plate bending theory is used to model the plate behavior and the Winkler foundation model is directly incorporated into the governing differential equation [11]. For this purpose, the concept of guided support edge (all the transverse shear force, the torsional moment and the average rotation of the normal to the mid-surface around the tangent of the plate boundary are equal to zero), was put forward. And accordingly, the original problem was decomposed into three subproblems. The first subproblem involves bending of the thick plate with four guided support edges and subjected to the transverse load, and the Navier solution was obtained. The other two subproblems involve bending of the thick plate with two opposite guided support edges (along *x*- or *y*-directions) and the other two edges subjected to distributed bending moments (along *y*- or *x*-directions). Under these circumstances, no transverse load is applied, and two Levy solutions were derived. Numerical examples show that the superposition solution exactly satisfies the boundary conditions of four free edges and has good computational accuracy. And furthermore, the size of the matrix of undetermined coefficients is just 1/36 of that in Henwood’s research, which implies an improved computational efficiency. This procedure was then employed for the analysis of thick Reissner plates with free edges and on Pasternak foundations [12]. When the boundary condition was extended (from the specific free edges) to the general case, the problem concerning the bending of a Reissner plate resting on the Pasternak foundation was also investigated by Sun et al. [13]. The basic equations were derived from Reissner theory, and modified to include the Pasternak foundation. This led to a system of three second order partial differential equations with respective to the deflection of the mid-surface, and its rotations about the axes. The series solution of the basic equations was arrived at with the general analytical Fourier series method. The general analytical (Fourier series) solution they obtained is of good convergence and good accuracy for not only the displacements and rotations, but also the higher order quantities such as bending moments. Furthermore, the prescribed types of boundary conditions (such as four clamped edges, four simply supported edges, and four free edges) are well fitted. Numerical investigations were performed for possible combinations of boundary conditions, loading conditions, thickness-to-length ratios, and parameters of thick plates and foundations. It was found that as the material parameter (Young’s modulus) of thick plates increases, the displacement decreases and is gradually stabilized, while the maximum tensile stress increases; as the thickness of thick plates increases, the displacement and maximum tensile stress decrease. These conclusions are constructive for the design, construction and maintain of concrete pavements of engineering structures.

In the last two decades, multiscale problems in science and engineering have attracted considerable attentions [14]. The spatial, material and/or physical parameters in a multiscale problem, can vary easily by orders of magnitude, resulting in the appearance of boundary layers or other forms of local discontinuities with large gradients in the solution region. For example, the convection-diffusion-reaction equation takes a very simple form of the second order linear differential equation. However, when the model parameters change within a wide spectrum, it has varied physical behaviors with the presence of the exponential regime and the propagation regime [15-18]. And it becomes one of the most studied subjects in multiscale methods. Compared to the convection-diffusion-reaction equation, the governing equations for elastic bending of Reissner plates on Pasternak foundations are much more complicated. Although when the load is positioned at the corner of the plate, the discontinuity of the shear force near the corner was observed to give rise to a significant localized change (large gradient) [10], the investigations were still confined within the specific structural analysis of engineering structures, and not extended to a thorough multiscale analysis for a wide spectrum of model parameters.

Nowadays, many multiscale methods have been developed for the multiscale problems, which include the stabilized finite element methods [19, 20], the bubble methods [21-24], the wavelet finite element methods [25-27], the meshless methods [28], the variational multiscale methods [29, 30], and so on. It is evident from their names that these multiscale methods are typically based on the traditional mesh-based or other discretization methods by making some suitable modifications, such as: use of stabilization terms, inclusion of different scale groups, decomposition of the solution into coarse and fine scale components, *etc*. Although these multiscale methods have been applied to a variety of boundary value problems, such as the convection-diffusion-reaction equation, we are still in the constant struggles with respect to the stability of the solution algorithms, proper selection of computational scales, robustness and effectiveness of the methods, balancing the span of scale groups and solution accuracy, high computational costs, low numerical accuracy for higher order derivatives of field variables, and so on.

It's worthy of notice that, a new type Fourier series method, namely the Fourier series method with supplementary terms, potentially represents a strategic shift from the existing framework towards better resolving multiscale problems. In this method, the solution function is expressed as a linear combination of a conventional Fourier series and some supplementary terms [31-49]. The supplementary terms, on the one hand, are purposely introduced to carry over the discontinuities potentially associated with the periodic extensions of the original solution function. As a result, the Fourier series actually corresponds to a periodic and sufficiently smooth residual function, and is hence uniformly convergent and termwise differentiable. On the other hand, the supplementary terms are only required to be sufficiently smooth over a compact interval (or domain), without regulating them to any particular forms. This implies that there is actually a large (theoretically, an infinite) number of possible choices for such basis functions in the process of implementing the Fourier series expansions. Specifically, if the supplementary terms are sought as the general solution of the differential equations, it can be expected that such general solutions shall be able to appropriately interpret the meaning of the differential equation, and hence better capture the spatial characteristics of the solution in separate directions. In view of the multiscale capability of this solution method, we rename it the Fourier series (based) multiscale method.

Therefore, in the fifth paper of the series of researches on Fourier series multiscale method, we will reinvestigate the bending problem of a Reissner plate on the Pasternak foundation from a more systematic and general view of point. Firstly, we generalize this issue to the multiscale analysis of a system of a fourth order linear differential equation (for transverse displacement of the plate) and a second order linear differential equation (for the stress function), where general boundary conditions and a wide spectrum of model parameters are prescribed. Secondly, the Fourier series multiscale method is developed for the analysis of the bending problem of a Reissner plate on the Pasternak foundation. As a routine task, the solution functions (for the transverse displacement and the stress function) are respectively decomposed into several constituents, such as the corner function, the two boundary functions and the internal function. The sum of the corner function and the internal function corresponds to the particular solutions. And the two boundary functions correspond to the general solutions which satisfy the homogeneous form of the governing differential equation. And then, with all the constituents selected respectively as polynomials, one-dimensional half-range Fourier series along the $y$ (or $x$)-direction, and two-dimensional half-range Fourier series, the Fourier series multiscale solution of the bending problem of a Reissner plate on the Pasternak foundation is derived. Accordingly, this paper begins with description of the problem. Detailed formulations related to the Fourier series multiscale method is then presented. Finally, convergence characteristics of the Fourier series multiscale solution are investigated with numerical examples, and the multiscale characteristics of the bending problem of a Reissner plate on the Pasternak foundation are demonstrated for a wide spectrum of model parameters.

## 2. Description of the problem

In this section, we make the following assumptions about the thick plate and the elastic foundations:

1. The thick plate is an elastic plate of uniform thickness with Young's modulus $E$, Poisson's ratio $\mu$, shear modulus $G$ and thickness $h$ over the domain $[0,a]\times[0,b]$, and its deformation and stress can be described by the deflection of mid-surface $w$ and average rotations formed by the normal to the mid-surface rotating around axis $x_2$ and axis $x_1$.

2. For the thick plate, the upper surface is subjected to a transverse load $q$, and the

lower surface only a reactive force $q_e$ of the foundation. A perfect contact exists between the plate and the foundation. In general, we have a biparametric representation [50] of the reactive force $q_e(x_1, x_2)$ of the foundation

$$q_e = kw - G_p \nabla^2 w\,, \tag{1}$$

where $k$ is the modulus of subgrade reaction of the foundation, $G_p$ is the shear modulus of the foundation and the Laplace differential operator

$$\nabla^2 = \frac{\partial^2}{\partial x_1^2} + \frac{\partial^2}{\partial x_2^2}\,. \tag{2}$$

Then, according to Reissner's thick plate theory [51], we can express the differential equation of equilibrium of transverse elastic bending of thick plate on biparametric foundations as

$$\left.\begin{aligned} \mathcal{L}_p w = (q - q_e) - \frac{(2-\mu)h^2}{10(1-\mu)} \nabla^2 (q - q_e) \\ \mathcal{L}_\psi \psi = 0 \end{aligned}\right\}, \tag{3}$$

where $\psi(x_1, x_2)$ is a stress function, and the differential operators

$$\mathcal{L}_p = D \nabla^2 \nabla^2\,, \tag{4}$$

$$\mathcal{L}_\psi = \nabla^2 - \frac{10}{h^2}\,, \tag{5}$$

and the flexural rigidity of the thick plate $D = Eh^3 / 12(1-\mu^2)$.

Then we have the physical equations

$$\begin{bmatrix} M_{x_1} \\ M_{x_2} \\ M_{x_1 x_2} \end{bmatrix} = D \begin{bmatrix} 1 & \mu & 0 \\ \mu & 1 & 0 \\ 0 & 0 & \frac{1-\mu}{2} \end{bmatrix} \begin{bmatrix} \frac{\partial \beta_{x_1}}{\partial x_1} \\ \frac{\partial \beta_{x_2}}{\partial x_2} \\ \frac{\partial \beta_{x_1}}{\partial x_2} + \frac{\partial \beta_{x_2}}{\partial x_1} \end{bmatrix}, \tag{6}$$

and

$$\begin{bmatrix} Q_{x_1} \\ Q_{x_2} \end{bmatrix} = \begin{bmatrix} C_s & 0 \\ 0 & C_s \end{bmatrix} \begin{bmatrix} \beta_{x_1} + \frac{\partial w}{\partial x_1} \\ \beta_{x_2} + \frac{\partial w}{\partial x_2} \end{bmatrix}, \tag{7}$$

where $M_{x_1}$, $M_{x_2}$ and $M_{x_1 x_2}$ are the bending moments and torsional moment of the thick plate, $Q_{x_1}$ and $Q_{x_2}$ are the transverse shear forces of the thick plate, and shear constant $C_s = 5Gh/6$.

In addition, the transverse shear forces $Q_{x_1}$ and $Q_{x_2}$ of the thick plate and the stress function $\psi$ satisfies the following equations

$$\left.\begin{aligned} Q_{x_1} = \frac{\partial \psi}{\partial x_2} - D \frac{\partial}{\partial x_1} (\nabla^2 w) \\ Q_{x_2} = -\frac{\partial \psi}{\partial x_1} - D \frac{\partial}{\partial x_2} (\nabla^2 w) \end{aligned}\right\}. \tag{8}$$

For example, there are three kinds of typical boundary conditions at the edge $x_1 = a$ of the plate:

1. Generalized clamped edge (C type edge for short)

$$w(a,x_2) = \bar{w}(a,x_2)\,,\;\; \beta_{x_1}(a,x_2) = \bar{\beta}_{x_1}(a,x_2)\,,\;\; \beta_{x_2}(a,x_2) = \bar{\beta}_{x_2}(a,x_2)\,, \tag{9}$$

2. Generalized simply supported edge (S type edge for short)

$$w(a,x_2) = \bar{w}(a,x_2)\,,\;\; M_{x_1}(a,x_2) = \bar{M}_{x_1}(a,x_2)\,,\;\; M_{x_1x_2}(a,x_2) = \bar{M}_{x_1x_2}(a,x_2)\,, \tag{10}$$

3. Generalized free edge (F type edge for short)

$$Q_{x_1}(a,x_2) = \bar{Q}_{x_1}(a,x_2)\,,\;\; M_{x_1}(a,x_2) = \bar{M}_{x_1}(a,x_2)\,,\;\; M_{x_1x_2}(a,x_2) = \bar{M}_{x_1x_2}(a,x_2)\,, \tag{11}$$

where $\bar{w}$, $\bar{\beta}_{x_1}$, $\bar{\beta}_{x_2}$, $\bar{M}_{x_1}$, $\bar{M}_{x_1x_2}$ and $\bar{Q}_{x_1}$ are the specified displacement, rotations, distributed bending moment, torsional moment and shear force, respectively.

## 3. The Fourier series multiscale solution

The elastic bending of thick plates on biparametric foundations can be formulated as a system of differential equations consisting of a fourth order ($2r_w = 4$) linear differential equations with constant coefficients of the transverse displacement of thick plate and a second order ($2r_\psi = 2$) linear differential equations with constant coefficients of the stress function. Moreover, the system of equations includes only the solution functions and their even order partial derivatives. Therefore, the transverse displacement of thick plate and the stress function will be sought in the forms of the two-dimensional half-range sine-sine series and the two-dimensional half-range cosine-cosine series, respectively [52].

### *3.1. The general solution of the transverse displacement*

We define the constants

$$D_h = D + \frac{(2-\mu)h^2}{10(1-\mu)}G_p\,,\;\; G_{ph} = G_p + \frac{(2-\mu)h^2}{10(1-\mu)}k\,,$$

and denote the differential operators by

$$\mathcal{L}_{pf} = D_h\nabla^2\nabla^2 - G_{ph}\nabla^2 + k\,, \tag{12}$$

$$\mathcal{L}_q = 1 - \frac{(2-\mu)h^2}{10(1-\mu)}\nabla^2\,. \tag{13}$$

We then substitute Eq. (1) in Eq. (3.a), and obtain

$$\mathcal{L}_{pf}\, w(x_1,x_2) = \mathcal{L}_q\, q(x_1,x_2)\,. \tag{14}$$

Suppose that the homogeneous solution of Eq. (14) is given as

$$p_{1n,Hw}(x_1,x_2) = \exp(\eta_{nw}x_1)\exp(\mathrm{i}\beta_n x_2)\,, \tag{15}$$

where $n$ is a positive integer, $\eta_{nw}$ is an undetermined constant, $\beta_n = n\pi/b$ and $\mathrm{i} = \sqrt{-1}$.

Substituting Eq. (15) into the homogeneous form of Eq. (14), we obtain the characteristics equation

$$D_h\eta_{nw}^4 - (2D_h\beta_n^2 + G_{ph})\eta_{nw}^2 + D_h\beta_n^4 + G_{ph}\beta_n^2 + k = 0\,. \tag{16}$$

It is easy to verify that Eq. (16) has the following four distinct real roots when $\Delta_h = G_{ph}^2 - 4D_h k > 0$,

$$\eta_{nw,1} = \alpha_{1n}\,,\ \ \eta_{nw,2} = -\alpha_{1n}\,,\ \ \eta_{nw,3} = \alpha_{2n}\,,\ \ \eta_{nw,4} = -\alpha_{2n}\,, \tag{17}$$

where $\alpha_{1n} = \sqrt{\beta_n^2 + \frac{1}{2D_h}\left[G_{ph} + \sqrt{G_{ph}^2 - 4D_h k}\right]}$,

$$\alpha_{2n} = \sqrt{\beta_n^2 + \frac{1}{2D_h}\left[G_{ph} - \sqrt{G_{ph}^2 - 4D_h k}\right]}\,.$$

Eq. (16) has two distinct double real roots when $\Delta_h = G_{ph}^2 - 4D_h k = 0$,

$$\eta_{nw,1} = \eta_{nw,2} = \alpha_{3n}\,,\ \ \eta_{nw,3} = \eta_{nw,4} = -\alpha_{3n}\,, \tag{18}$$

where $\alpha_{3n} = \sqrt{\beta_n^2 + \frac{G_{ph}}{2D_h}}$ .

Eq. (16) has four distinct complex roots when $\Delta_h = G_{ph}^2 - 4D_h k < 0$

$$\eta_{nw,1,2,3,4} = \pm\alpha_{5n} \pm \mathrm{i}\alpha_{6n}\,, \tag{19}$$

where $\alpha_{5n} > 0$, $\alpha_{6n} > 0$ and $\alpha_{5n} + \mathrm{i}\alpha_{6n} = \sqrt{\beta_n^2 + \frac{1}{2D_h}\left[G_{ph} + \mathrm{i}\sqrt{4D_h k - G_{ph}^2}\right]}$.

With the expressions of $p_{1nl,w}(x_1)$, $l = 1, 2, 3, 4$, presented in Table 1, we obtain the expression of $p_{1nl,Hw}(x_1, x_2)$ as follows

$$p_{1nl,Hw}(x_1, x_2) = p_{1nl,w}(x_1)\sin(\beta_n x_2)\,,\ \ l = 1, 2, 3, 4\,. \tag{20}$$

Table 1: Expressions for $p_{1nl,w}(x_1)$, $l = 1, 2, 3, 4$.

| | $\Delta_h > 0$ | $\Delta_h = 0$ | $\Delta_h < 0$ |
|---|---|---|---|
| $p_{1n1,w}(x_1)$ | $\frac{\sinh(\alpha_{1n}x_1)}{\sinh(\alpha_{1n}a)}$ | $\frac{\sinh(\alpha_{3n}x_1)}{\sinh(\alpha_{3n}a)}$ | $\frac{\sinh(\alpha_{5n}x_1)\sin(\alpha_{6n}x_1)}{\sinh(\alpha_{5n}a)\sin(\alpha_{6n}a)}$ |
| $p_{1n2,w}(x_1)$ | $\frac{\sinh[\alpha_{1n}(a-x_1)]}{\sinh(\alpha_{1n}a)}$ | $\frac{x_1\sinh(\alpha_{3n}x_1)}{a\sinh(\alpha_{3n}a)}$ | $\frac{\sinh(\alpha_{5n}x_1)\sin[\alpha_{6n}(a-x_1)]}{\sinh(\alpha_{5n}a)\sin(\alpha_{6n}a)}$ |
| $p_{1n3,w}(x_1)$ | $\frac{\sinh(\alpha_{2n}x_1)}{\sinh(\alpha_{2n}a)}$ | $\frac{\sinh[\alpha_{3n}(a-x_1)]}{\sinh(\alpha_{3n}a)}$ | $\frac{\sinh[\alpha_{5n}(a-x_1)]\sin(\alpha_{6n}x_1)}{\sinh(\alpha_{5n}a)\sin(\alpha_{6n}a)}$ |
| $p_{1n4,w}(x_1)$ | $\frac{\sinh[\alpha_{2n}(a-x_1)]}{\sinh(\alpha_{2n}a)}$ | $\frac{(a-x_1)\sinh[\alpha_{3n}(a-x_1)]}{a\sinh(\alpha_{3n}a)}$ | $\frac{\sinh[\alpha_{5n}(a-x_1)]\sin[\alpha_{6n}(a-x_1)]}{\sinh(\alpha_{5n}a)\sin(\alpha_{6n}a)}$ |

Selecting the following vector of functions

$$\mathbf{p}_{1n,Hw}^{\mathrm{T}}(x_1, x_2) = [\,p_{1n1,Hw}(x_1, x_2)\quad p_{1n2,Hw}(x_1, x_2)\quad p_{1n3,Hw}(x_1, x_2)\quad p_{1n4,Hw}(x_1, x_2)]\,, \tag{21}$$

then we can express the general solution of Eq. (14) as

$$w_1(x_1, x_2) = \mathbf{\Phi}_{1,w}^{\mathrm{T}}(x_1, x_2)\cdot\mathbf{q}_{1,w}\,, \tag{22}$$

where $\mathbf{\Phi}_{1,w}^{\mathrm{T}}(x_1, x_2)$, a vector of functions, and $\mathbf{q}_{1,w}$, an undetermined coefficient vector, are defined as in (92)-(100) in [52]. However, the relevant expressions need to be properly modified with the aid of Theorem 7 in [53].

*3.2. The general solution of the stress function*

Suppose that Eq. (3.b) has the following homogeneous solution

$$p_{1n,H\psi}(x_1,x_2)=\exp(\eta_{n\psi}x_1)\exp(\mathrm{i}\beta_n x_2)\,, \tag{23}$$

where $n$ is a positive integer, $\eta_{n\psi}$ is an undetermined constant, $\beta_n = n\pi/b$ and $\mathrm{i}=\sqrt{-1}$.

Substituting Eq. (23) in Eq. (3.b), we obtain the characteristics equation

$$\eta_{n\psi}^2-\beta_n^2-\frac{10}{h^2}=0\,. \tag{24}$$

It is easy to find that Eq. (24) has two real distinct roots

$$\eta_{n\psi,1}=\alpha_{7n}\,,\quad \eta_{n\psi,2}=-\alpha_{7n}\,, \tag{25}$$

where $\alpha_{7n}=\sqrt{\beta_n^2+\dfrac{10}{h^2}}$ .

We define the one-dimensional basis functions

$$p_{1n1,\psi}(x_1)=\frac{\sinh(\alpha_{7n}x_1)}{\sinh(\alpha_{7n}a)}\,, \tag{26}$$

$$p_{1n2,\psi}(x_1)=\frac{\sinh[(\alpha_{7n}(a-x_1)]}{\sinh(\alpha_{7n}a)}\,, \tag{27}$$

and then $p_{1nl,H\psi}(x_1,x_2)$ has the following expression

$$p_{1nl,H\psi}(x_1,x_2)=p_{1nl,\psi}(x_1)\cos(\beta_n x_2)\,,\quad l=1,2\,. \tag{28}$$

Selecting the following vector of functions

$$\mathbf{p}_{1n,H\psi}^{\mathrm{T}}(x_1,x_2)=[p_{1n1,H\psi}(x_1,x_2)\quad p_{1n2,H\psi}(x_1,x_2)]\,, \tag{29}$$

then we can express the general solution of Eq. (3.b) as

$$\psi_1(x_1,x_2)=\mathbf{\Phi}_{1,\psi}^{\mathrm{T}}(x_1,x_2)\cdot\mathbf{q}_{1,\psi}\,, \tag{30}$$

where $\mathbf{\Phi}_{1,\psi}^{\mathrm{T}}(x_1,x_2)$, a vector of functions, and $\mathbf{q}_{1,\psi}$, an undetermined coefficient vector, are also defined as in (92)-(100) in [52]. However, the relevant expressions need to be properly modified with the aid of Theorem 8 in [53].

*3.3. Expression of the Fourier series multiscale solution*

According to the Fourier series multiscale method [52] of fourth order linear differential equation with constant coefficients, we expand the transverse displacement function $w(x_1,x_2)$ in a composite half-range sine-sine series over the domain $[0,a]\times[0,b]$

$$\begin{aligned} w(x_1,x_2)&=w_0(x_1,x_2)+w_1(x_1,x_2)+w_2(x_1,x_2)+w_3(x_1,x_2)\\ &=\mathbf{\Phi}_{0,w}^{\mathrm{T}}(x_1,x_2)\cdot\mathbf{q}_{0,w}+\mathbf{\Phi}_{1,w}^{\mathrm{T}}(x_1,x_2)\cdot\mathbf{q}_{1,w}+\mathbf{\Phi}_{2,w}^{\mathrm{T}}(x_1,x_2)\cdot\mathbf{q}_{2,w}+\mathbf{\Phi}_{3,w}^{\mathrm{T}}(x_1,x_2)\cdot\mathbf{q}_{3,w}\,. \end{aligned} \tag{31}$$

If we define the vector of basis functions

$$\mathbf{\Phi}_{w}^{\mathrm{T}}(x_1,x_2)=[\mathbf{\Phi}_{0,w}^{\mathrm{T}}(x_1,x_2)\quad \mathbf{\Phi}_{1,w}^{\mathrm{T}}(x_1,x_2)\quad \mathbf{\Phi}_{2,w}^{\mathrm{T}}(x_1,x_2)\quad \mathbf{\Phi}_{3,w}^{\mathrm{T}}(x_1,x_2)]\,, \tag{32}$$

and the undetermined constant vector

$$\mathbf{q}_{w}^{\mathrm{T}}=[\mathbf{q}_{0,w}^{\mathrm{T}}\quad \mathbf{q}_{1,w}^{\mathrm{T}}\quad \mathbf{q}_{2,w}^{\mathrm{T}}\quad \mathbf{q}_{3,w}^{\mathrm{T}}]\,, \tag{33}$$

then the displacement can be rewritten as

$$w(x_1,x_2)=\mathbf{\Phi}_{w}^{\mathrm{T}}(x_1,x_2)\cdot\mathbf{q}_{w}\,. \tag{34}$$

Similarly, according to the Fourier series multiscale method [52] of second order linear differential equation with constant coefficients, we expand the stress function $\psi(x_1, x_2)$ in a composite half-range cosine-cosine series over the domain $[0,a]\times[0,b]$

$$\begin{aligned}\psi(x_1,x_2) &= \psi_0(x_1,x_2)+\psi_1(x_1,x_2)+\psi_2(x_1,x_2)\\ &= \mathbf{\Phi}_{0,\psi}^{\mathrm{T}}(x_1,x_2)\cdot\mathbf{q}_{0,\psi}+\mathbf{\Phi}_{1,\psi}^{\mathrm{T}}(x_1,x_2)\cdot\mathbf{q}_{1,\psi}+\mathbf{\Phi}_{2,\psi}^{\mathrm{T}}(x_1,x_2)\cdot\mathbf{q}_{2,\psi}\,.\end{aligned}\tag{35}$$

Since Eq. (3.b) is a homogeneous equation, we have

$$\mathbf{q}_{0,\psi}=\mathbf{0}\,.\tag{36}$$

If we define the vector of basis functions

$$\mathbf{\Phi}_{\psi}^{\mathrm{T}}(x_1,x_2)=[\mathbf{\Phi}_{1,\psi}^{\mathrm{T}}(x_1,x_2)\quad \mathbf{\Phi}_{2,\psi}^{\mathrm{T}}(x_1,x_2)]\,,\tag{37}$$

and the undetermined constant vector

$$\mathbf{q}_{\psi}^{\mathrm{T}}=[\mathbf{q}_{1,\psi}^{\mathrm{T}}\quad \mathbf{q}_{2,\psi}^{\mathrm{T}}]\,,\tag{38}$$

then the stress function can be rewritten as

$$\psi(x_1,x_2)=\mathbf{\Phi}_{\psi}^{\mathrm{T}}(x_1,x_2)\cdot\mathbf{q}_{\psi}\,.\tag{39}$$

By combining Eqs. (34) and (39), we then obtain the Fourier series multiscale solution for elastic bending of thick plates on biparametric foundations

$$\begin{bmatrix} w\\ \psi\end{bmatrix}=\mathbf{\Phi}\cdot\mathbf{q}\,,\tag{40}$$

where the matrix of the basis functions

$$\mathbf{\Phi}=\begin{bmatrix}\mathbf{\Phi}_{w}^{\mathrm{T}} & 0\\ 0 & \mathbf{\Phi}_{\psi}^{\mathrm{T}}\end{bmatrix},\tag{41}$$

and the vector of the undetermined constants

$$\mathbf{q}^{\mathrm{T}}=[\mathbf{q}_{w}^{\mathrm{T}}\quad \mathbf{q}_{\psi}^{\mathrm{T}}]\,.\tag{42}$$

### *3.4. Equivalent transformation of the solution*

Substituting the Fourier series multiscale solution in Eqs. (7) and (8), we obtain

$$\begin{bmatrix} w\\ \beta_{x_1}\\ \beta_{x_2}\end{bmatrix}=\mathbf{\Gamma}\cdot\mathbf{q}\,,\tag{43}$$

where the transformation matrix

$$\mathbf{\Gamma}=\begin{bmatrix}\mathbf{\Phi}_{w}^{\mathrm{T}} & \mathbf{0}\\ -\dfrac{D}{C_s}\mathbf{\Phi}_{w}^{(3,0)\mathrm{T}}-\dfrac{D}{C_s}\mathbf{\Phi}_{w}^{(1,2)\mathrm{T}}-\mathbf{\Phi}_{w}^{(1,0)\mathrm{T}} & \dfrac{1}{C_s}\mathbf{\Phi}_{\psi}^{(0,1)\mathrm{T}}\\ -\dfrac{D}{C_s}\mathbf{\Phi}_{w}^{(2,1)\mathrm{T}}-\dfrac{D}{C_s}\mathbf{\Phi}_{w}^{(0,3)\mathrm{T}}-\mathbf{\Phi}_{w}^{(0,1)\mathrm{T}} & -\dfrac{1}{C_s}\mathbf{\Phi}_{\psi}^{(1,0)\mathrm{T}}\end{bmatrix}.\tag{44}$$

And specifically for the edge $x_1=a$, we obtain

$$\begin{bmatrix} w(a,x_2)\\ \beta_{x_1}(a,x_2)\\ \beta_{x_2}(a,x_2)\end{bmatrix}=\mathbf{\Gamma}(a,x_2)\cdot\mathbf{q}\,.\tag{45}$$

Expand $w(a,x_2)$, $\beta_{x_1}(a,x_2)$, $\beta_{x_2}(a,x_2)$ and the corresponding elements in the first, second and third rows of the matrix $\mathbf{\Gamma}(a,x_2)$ respectively in half-range sine series, half-range sine series and half-range cosine series over the interval $(0,b)$, and let $N$ be the number of truncated terms. Then by comparing the Fourier coefficients of like terms in the left-hand side of Eq. (45) with that of the corresponding terms in the right-hand side successively, we then obtain

$$\begin{bmatrix} \mathbf{q}_{w,1a} \\ \mathbf{q}_{\beta x_1,1a} \\ \mathbf{q}_{\beta x_2,1a} \end{bmatrix} = \mathbf{R}_{pf,1a} \cdot \mathbf{q}, \tag{46}$$

where $\mathbf{q}_{w,1a}$, $\mathbf{q}_{\beta x_1,1a}$ and $\mathbf{q}_{\beta x_2,1a}$ are Fourier coefficient vectors corresponding to $w(a,x_2)$, $\beta_{x_1}(a,x_2)$ and $\beta_{x_2}(a,x_2)$, or might be called the Fourier coefficient subvectors of boundary conditions corresponding to the edge $x_1 = a$, and $\mathbf{R}_{pf,1a}$ is the matrix of Fourier coefficients corresponding to $\mathbf{\Gamma}(a,x_2)$.

Similarly, for the edges $x_1 = 0$, $x_2 = b$ and $x_2 = 0$, we obtain another three systems of equations. Combining them with Eq. (46), we have

$$\mathbf{q}_b = \mathbf{R}_{pf} \cdot \mathbf{q}, \tag{47}$$

where the vector of boundary Fourier coefficients

$$\mathbf{q}_b = \begin{bmatrix} \mathbf{q}_{b,1a} \\ \mathbf{q}_{b,10} \\ \mathbf{q}_{b,2b} \\ \mathbf{q}_{b,20} \end{bmatrix}, \tag{48}$$

with the subvectors of boundary Fourier coefficients

$$\mathbf{q}_{b,1a} = \begin{bmatrix} \mathbf{q}_{w,1a} \\ \mathbf{q}_{\beta x_1,1a} \\ \mathbf{q}_{\beta x_2,1a} \end{bmatrix}, \quad \mathbf{q}_{b,10} = \begin{bmatrix} \mathbf{q}_{w,10} \\ \mathbf{q}_{\beta x_1,10} \\ \mathbf{q}_{\beta x_2,10} \end{bmatrix}, \quad \mathbf{q}_{b,2b} = \begin{bmatrix} \mathbf{q}_{w,2b} \\ \mathbf{q}_{\beta x_1,2b} \\ \mathbf{q}_{\beta x_2,2b} \end{bmatrix}, \quad \mathbf{q}_{b,20} = \begin{bmatrix} \mathbf{q}_{w,20} \\ \mathbf{q}_{\beta x_1,20} \\ \mathbf{q}_{\beta x_2,20} \end{bmatrix}, \tag{49}$$

and the matrix

$$\mathbf{R}_{pf} = \begin{bmatrix} \mathbf{R}_{pf,1a} \\ \mathbf{R}_{pf,10} \\ \mathbf{R}_{pf,2b} \\ \mathbf{R}_{pf,20} \end{bmatrix}. \tag{50}$$

We define the vectors of undetermined constants

$$\mathbf{q}_{03}^{\mathrm{T}} = [\mathbf{q}_{0,w}^{\mathrm{T}} \quad \mathbf{q}_{3,w}^{\mathrm{T}}], \tag{51}$$

$$\mathbf{q}_{12}^{\mathrm{T}} = [\mathbf{q}_{1,w}^{\mathrm{T}} \quad \mathbf{q}_{2,w}^{\mathrm{T}} \quad \mathbf{q}_{1,\psi}^{\mathrm{T}} \quad \mathbf{q}_{2,\psi}^{\mathrm{T}}], \tag{52}$$

and write the matrix $\mathbf{R}_{pf}$ in block form

$$\mathbf{R}_{pf} = [\mathbf{R}_{pf,03} \quad \mathbf{R}_{pf,12}], \tag{53}$$

then Eq. (47) can be rewritten as

$$\mathbf{q}_b = \mathbf{R}_{pf,03} \cdot \mathbf{q}_{03} + \mathbf{R}_{pf,12} \cdot \mathbf{q}_{12}. \tag{54}$$

Therefore, we have

$$\mathbf{q}_{12} = -\mathbf{R}_{pf,12}^{-1}\mathbf{R}_{pf,03} \cdot \mathbf{q}_{03} + \mathbf{R}_{pf,12}^{-1} \cdot \mathbf{q}_b. \tag{55}$$

Substituting it in Eq. (40), the Fourier series multiscale solution for elastic bending of thick plates on biparametric foundations can be rewritten as

$$\begin{bmatrix} w \\ \psi \end{bmatrix} = \mathbf{\Phi}_R \cdot \mathbf{q}_R \,, \tag{56}$$

where the modified matrix of basis functions

$$\mathbf{\Phi}_R = \begin{bmatrix} \mathbf{\Phi}_{0,w}^{\mathrm{T}} & \mathbf{\Phi}_{3,w}^{\mathrm{T}} & \mathbf{\Phi}_{1,w}^{\mathrm{T}} & \mathbf{\Phi}_{2,w}^{\mathrm{T}} & \mathbf{0} & \mathbf{0} \\ \mathbf{0} & \mathbf{0} & \mathbf{0} & \mathbf{0} & \mathbf{\Phi}_{1,\psi}^{\mathrm{T}} & \mathbf{\Phi}_{2,\psi}^{\mathrm{T}} \end{bmatrix} \begin{bmatrix} \mathbf{I} & \mathbf{0} \\ -\mathbf{R}_{pf,12}^{-1}\mathbf{R}_{pf,03} & \mathbf{R}_{pf,12}^{-1} \end{bmatrix}, \tag{57}$$

and the modified vector of undetermined constants

$$\mathbf{q}_R^{\mathrm{T}} = [\mathbf{q}_{03}^{\mathrm{T}} \quad \mathbf{q}_b^{\mathrm{T}}]\,. \tag{58}$$

### *3.5. Expression of stress resultants*

We write the modified matrix of basis functions $\mathbf{\Phi}_R$ in block form

$$\mathbf{\Phi}_R = \begin{bmatrix} \mathbf{\Phi}_{Rw} \\ \mathbf{\Phi}_{R\psi} \end{bmatrix}. \tag{59}$$

Then it is easy to obtain the expression of transverse displacement and rotations of the thick plate

$$\begin{bmatrix} w \\ \beta_{x_1} \\ \beta_{x_2} \end{bmatrix} = \mathbf{\Gamma}_{Rd} \cdot \mathbf{q}_R \,, \tag{60}$$

where the matrix

$$\mathbf{\Gamma}_{Rd} = \begin{bmatrix} \mathbf{\Gamma}_{R,w} \\ \mathbf{\Gamma}_{R,\beta_{x1}} \\ \mathbf{\Gamma}_{R,\beta_{x2}} \end{bmatrix} = \begin{bmatrix} \mathbf{\Phi}_{Rw} \\ -\frac{D}{C_s}\mathbf{\Phi}_{Rw}^{(3,0)} - \frac{D}{C_s}\mathbf{\Phi}_{Rw}^{(1,2)} - \mathbf{\Phi}_{Rw}^{(1,0)} + \frac{1}{C_s}\mathbf{\Phi}_{R\psi}^{(0,1)} \\ -\frac{D}{C_s}\mathbf{\Phi}_{Rw}^{(2,1)} - \frac{D}{C_s}\mathbf{\Phi}_{Rw}^{(0,3)} - \mathbf{\Phi}_{Rw}^{(0,1)} - \frac{1}{C_s}\mathbf{\Phi}_{R\psi}^{(1,0)} \end{bmatrix}. \tag{61}$$

Further, we obtain the expression of bending moments and torsional moment of the thick plate

$$\begin{bmatrix} M_{x_1} \\ M_{x_2} \\ M_{x_1x_2} \end{bmatrix} = \mathbf{\Gamma}_{RM} \cdot \mathbf{q}_R \,, \tag{62}$$

where the matrix

$$\mathbf{\Gamma}_{RM} = \begin{bmatrix} \mathbf{\Gamma}_{R,M_{x1}} \\ \mathbf{\Gamma}_{R,M_{x2}} \\ \mathbf{\Gamma}_{R,M_{x1x2}} \end{bmatrix}$$

$$
= \begin{bmatrix} -\frac{D^2}{C_s}\mathbf{\Phi}_{Rw}^{(4,0)} - \frac{D^2(1+\mu)}{C_s}\mathbf{\Phi}_{Rw}^{(2,2)} - D\mathbf{\Phi}_{Rw}^{(2,0)} - \frac{\mu D^2}{C_s}\mathbf{\Phi}_{Rw}^{(0,4)} - \mu D\mathbf{\Phi}_{Rw}^{(0,2)} + \frac{D(1-\mu)}{C_s}\mathbf{\Phi}_{R\psi}^{(1,1)} \\ -\frac{\mu D^2}{C_s}\mathbf{\Phi}_{Rw}^{(4,0)} - \frac{D^2(1+\mu)}{C_s}\mathbf{\Phi}_{Rw}^{(2,2)} - \mu D\mathbf{\Phi}_{Rw}^{(2,0)} - \frac{D^2}{C_s}\mathbf{\Phi}_{Rw}^{(0,4)} - D\mathbf{\Phi}_{Rw}^{(0,2)} - \frac{D(1-\mu)}{C_s}\mathbf{\Phi}_{R\psi}^{(1,1)} \\ -\frac{D^2(1-\mu)}{C_s}\mathbf{\Phi}_{Rw}^{(3,1)} - \frac{D^2(1-\mu)}{C_s}\mathbf{\Phi}_{Rw}^{(1,3)} - D(1-\mu)\mathbf{\Phi}_{Rw}^{(1,1)} - \frac{D(1-\mu)}{2C_s}\mathbf{\Phi}_{R\psi}^{(2,0)} + \frac{D(1-\mu)}{2C_s}\mathbf{\Phi}_{R\psi}^{(0,2)} \end{bmatrix}. \tag{63}
$$

Meanwhile, we obtain the expression of shear forces of the thick plate

$$
\begin{bmatrix} Q_{x_1} \\ Q_{x_2} \end{bmatrix} = \mathbf{\Gamma}_{RQ} \cdot \mathbf{q}_R, \tag{64}
$$

where the matrix

$$
\mathbf{\Gamma}_{RQ} = \begin{bmatrix} \mathbf{\Gamma}_{R,Q_{x1}} \\ \mathbf{\Gamma}_{R,Q_{x2}} \end{bmatrix} = \begin{bmatrix} -D\mathbf{\Phi}_{Rw}^{(3,0)} - D\mathbf{\Phi}_{Rw}^{(1,2)} + \mathbf{\Phi}_{R\psi}^{(0,1)} \\ -D\mathbf{\Phi}_{Rw}^{(2,1)} - D\mathbf{\Phi}_{Rw}^{(0,3)} - \mathbf{\Phi}_{R\psi}^{(1,0)} \end{bmatrix}. \tag{65}
$$

## 4. Problem solving

For the elastic system consisting of thick Reissner plate and biparametric foundation, we can formulate the energy of the system as [3, 8]

$$
\Pi = \Pi_p + \Pi_f + \Pi_q + \Pi_\sigma, \tag{66}
$$

where $\Pi_p$ is the elastic potential energy of thick Reissner plate, $\Pi_f$ is the elastic potential energy of biparametric foundation, $\Pi_q$ is the potential energy of transverse load distributed over the thick plate, and $\Pi_\sigma$ is the total potential energy of bending moment, torsional moment and transverse shear force distributed on the stress boundaries of the thick plate.

Specifically, the elastic potential energy of thick Reissner plate is

$$
\Pi_p = \frac{1}{2}\int_0^a\int_0^b [d_{11}M_{x_1}^2 + 2d_{12}M_{x_1}M_{x_2} + d_{22}M_{x_2}^2 + d_{33}M_{x_1x_2}^2 + d_{44}Q_{x_1}^2 + d_{55}Q_{x_2}^2]dx_1dx_2, \tag{67}
$$

where the constants

$$
d_{11} = d_{22} = \frac{1}{D(1-\mu^2)},\ d_{12} = -\frac{\mu}{D(1-\mu^2)},\ d_{33} = \frac{2}{D(1-\mu)},\ d_{44} = d_{55} = \frac{1}{C_s}.
$$

The elastic potential energy of biparametric foundations is

$$
\Pi_p = \frac{1}{2}\int_0^a\int_0^b [kw^2 + G_p(\frac{\partial w}{\partial x_1})^2 + G_p(\frac{\partial w}{\partial x_2})^2]dx_1dx_2. \tag{68}
$$

The potential energy of transverse load distributed over the thick plate is

$$
\Pi_p = -\int_0^a\int_0^b qwdx_1dx_2. \tag{69}
$$

If all edges of the thick plate are stress boundaries and there are bending moment, torsional moment and transverse shear force distributed along the edges, the corresponding potential energy is

$$
\Pi_\sigma = -\int_0^b [\overline{Q}_{x_1}(a,x_2)w(a,x_2) + \overline{M}_{x_1}(a,x_2)\beta_{x_1}(a,x_2) + \overline{M}_{x_1x_2}(a,x_2)\beta_{x_2}(a,x_2)]dx_2
$$

$$+\int_0^b [\overline{Q}_{x_1}(0,x_2)w(0,x_2)+\overline{M}_{x_1}(0,x_2)\beta_{x_1}(0,x_2)+\overline{M}_{x_1x_2}(0,x_2)\beta_{x_2}(0,x_2)]dx_2$$
$$-\int_0^a [\overline{Q}_{x_2}(x_1,b)w(x_1,b)+\overline{M}_{x_2}(x_1,b)\beta_{x_2}(x_1,b)+\overline{M}_{x_1x_2}(x_1,b)\beta_{x_1}(x_1,b)]dx_1$$
$$+\int_0^a [\overline{Q}_{x_2}(x_1,0)w(x_1,0)+\overline{M}_{x_2}(x_1,0)\beta_{x_2}(x_1,0)+\overline{M}_{x_1x_2}(x_1,0)\beta_{x_1}(x_1,0)]dx_1 . \tag{70}$$

Substituting Eqs. (67)-(70) in Eq. (66), we rewrite Eq. (66) as

$$\Pi = \frac{1}{2}\mathbf{q}_R^{\mathrm{T}}\mathbf{K}_{pf}\mathbf{q}_R - \mathbf{q}_R^{\mathrm{T}}\mathbf{Q}_{pf} , \tag{71}$$

where the stiffness matrix

$$\mathbf{K}_{pf} = \int_0^a\int_0^b [d_{11}\mathbf{\Gamma}^{\mathrm{T}}_{R,M_{x1}}\mathbf{\Gamma}_{R,M_{x1}} + d_{12}\mathbf{\Gamma}^{\mathrm{T}}_{R,M_{x1}}\mathbf{\Gamma}_{R,M_{x2}} + d_{12}\mathbf{\Gamma}^{\mathrm{T}}_{R,M_{x2}}\mathbf{\Gamma}_{R,M_{x1}} + d_{22}\mathbf{\Gamma}^{\mathrm{T}}_{R,M_{x2}}\mathbf{\Gamma}_{R,M_{x2}}$$
$$+d_{33}\mathbf{\Gamma}^{\mathrm{T}}_{R,M_{x1x2}}\mathbf{\Gamma}_{R,M_{x1x2}} + d_{44}\mathbf{\Gamma}^{\mathrm{T}}_{R,Q_{x1}}\mathbf{\Gamma}_{R,Q_{x1}} + d_{55}\mathbf{\Gamma}^{\mathrm{T}}_{R,Q_{x2}}\mathbf{\Gamma}_{R,Q_{x2}} + k\mathbf{\Gamma}^{\mathrm{T}}_{R,w}\mathbf{\Gamma}_{R,w}$$
$$+G_p(\frac{1}{C_s}\mathbf{\Gamma}_{R,Q_{x1}} - \mathbf{\Gamma}_{R,\beta_{x1}})^{\mathrm{T}}(\frac{1}{C_s}\mathbf{\Gamma}_{R,Q_{x1}} - \mathbf{\Gamma}_{R,\beta_{x1}})$$
$$+G_p(\frac{1}{C_s}\mathbf{\Gamma}_{R,Q_{x2}} - \mathbf{\Gamma}_{R,\beta_{x2}})^{\mathrm{T}}(\frac{1}{C_s}\mathbf{\Gamma}_{R,Q_{x2}} - \mathbf{\Gamma}_{R,\beta_{x2}})]dx_1dx_2 , \tag{72}$$

and the matrix of equivalent force

$$\mathbf{Q}_{pf} = \int_0^a\int_0^b q\mathbf{\Gamma}^{\mathrm{T}}_{R,w}$$
$$+\int_0^b [\overline{Q}_{x_1}(a,x_2)\mathbf{\Gamma}^{\mathrm{T}}_{R,w}(a,x_2)+\overline{M}_{x_1}(a,x_2)\mathbf{\Gamma}^{\mathrm{T}}_{R,\beta_{x1}}(a,x_2)+\overline{M}_{x_1x_2}(a,x_2)\mathbf{\Gamma}^{\mathrm{T}}_{R,\beta_{x2}}(a,x_2)]dx_2$$
$$-\int_0^b [\overline{Q}_{x_1}(0,x_2)\mathbf{\Gamma}^{\mathrm{T}}_{R,w}(0,x_2)+\overline{M}_{x_1}(0,x_2)\mathbf{\Gamma}^{\mathrm{T}}_{R,\beta_{x1}}(0,x_2)+\overline{M}_{x_1x_2}(0,x_2)\mathbf{\Gamma}^{\mathrm{T}}_{R,\beta_{x2}}(0,x_2)]dx_2$$
$$+\int_0^a [\overline{Q}_{x_2}(x_1,b)\mathbf{\Gamma}^{\mathrm{T}}_{R,w}(x_1,b)+\overline{M}_{x_2}(x_1,b)\mathbf{\Gamma}^{\mathrm{T}}_{R,\beta_{x2}}(x_1,b)+\overline{M}_{x_1x_2}(x_1,b)\mathbf{\Gamma}^{\mathrm{T}}_{R,\beta_{x1}}(x_1,b)]dx_1$$
$$-\int_0^a [\overline{Q}_{x_2}(x_1,0)\mathbf{\Gamma}^{\mathrm{T}}_{R,w}(x_1,0)+\overline{M}_{x_2}(x_1,0)\mathbf{\Gamma}^{\mathrm{T}}_{R,\beta_{x2}}(x_1,0)+\overline{M}_{x_1x_2}(x_1,0)\mathbf{\Gamma}^{\mathrm{T}}_{R,\beta_{x1}}(x_1,0)]dx_1 . \tag{73}$$

The corresponding stationary condition finally becomes

$$\mathbf{K}_{pf}\mathbf{q}_R = \mathbf{Q}_{pf} . \tag{74}$$

## 5. Numerical examples

In this section, the inverse validation method is adopted in the numerical examples of elastic bending of thick plates on the biparametric foundations. We select the following linear combination of homogeneous solutions $p_{1nl,Hw}(x_1,x_2)$ , $l=1,2,3,4$, in section 3.1 for the reference solution $w_{ref}(x_1,x_2)$

$$w_{ref}(x_1,x_2) = \sum_{l=1}^{4} 10^{-3} a \cdot p_{1nl,Hw}(x_1,x_2) + \frac{D}{ka^3} , \tag{75}$$

where the parameter $\beta_n = n\pi/b$ is substituted with $\beta_{ref} = \pi/2b$.

Meanwhile, the transverse load is the uniform load over the domain $[0,a]\times[0,b]$

$$q_{ref}(x_1,x_2) = Da^{-3} . \tag{76}$$

And in the setting of Fourier series multiscale computational scheme, the derivation method for discrete system is the variational method.

*5.1. Convergence characteristics*

As shown in Table 2, we have specially designed four sets of comparative numerical experiments for a detailed analysis of the influences of some key factors on the convergence characteristics and approximation accuracy of the Fourier series multiscale solution of the elastic bending of thick plates on biparametric foundations. These factors include computational parameters, thickness-to-length ratio, length-to-width ratio and boundary conditions. Based on the reference computational scheme, we perform the first set of numerical experiment by adjusting the computational parameters, $k_r = ka^4/D$ and $G_{pr} = G_p a^2/D$, of the foundation from 1and 1 to 100 and 10, $10^4$ and 100, and $10^6$ and 2000 successively. Based on the reference computational scheme, we perform the second set of numerical experiment by adjusting thickness-to-length ratio of the plate from 0.1 to 0.01, 0.2 and 0.4 successively. Based on the reference computational scheme, we perform the third set of numerical experiment by adjusting length-to-width ratio of the plate from 1.0 to 0.67, 0.50, 1.25 and 2.0 successively. Based on the reference computational scheme, we perform the fourth set of numerical experiment by adjusting boundary conditions from generalized clamped boundary condition (CCCC) to generalized simply supported boundary condition (SSSS) and generalized free boundary condition (FFFF) successively.

Table 2: Computational schemes for elastic bending of thick plates on biparametric foundations.

| Numerical experiment | No. | Boundary condition | Computational parameter ($k_r$, $G_{pr}$) | Length-to -width ratio $a/b$ | Thickness-to length ratio $h/a$ |
|---|---|---|---|---|---|
| 1 | a | CCCC | (1, 1) | 1.0 | 0.1 |
| | b | | (100, 10) | | |
| | c | | ($10^4$, 100) | | |
| | d | | ($10^6$, 2000) | | |
| 2 | a | CCCC | (1, 1) | 1.0 | 0.1 |
| | b | | | | 0.01 |
| | c | | | | 0.2 |
| | d | | | | 0.4 |
| 3 | a | CCCC | (1, 1) | 1.0 | 0.1 |
| | b | | | 0.67 | |
| | c | | | 0.50 | |
| | d | | | 1.25 | |
| | e | | | 2.0 | |
| | a | CCCC | | | |
| 4 | b | SSSS | (1, 1) | 1.0 | 0.1 |
| | c | FFFF | | | |

As to the reference computational scheme and four sets of comparative numerical experiments as given in Table 2, we take the first 2, 3, 5, 10, 15 and 20 terms successively in the composite Fourier series of the Fourier series multiscale solution and obtain the overall approximation errors, internal approximation errors, boundary approximation errors and corner approximation errors of the deflection of mid-surface $w$, average rotations $\beta_{x_1}$ and $\beta_{x_2}$, and bending moments $M_{x_1}$ and $M_{x_2}$ of the thick plate. Some results are presented in Figures 1-5. Convergence characteristics of the Fourier series multiscale solution of elastic bending of thick plates on biparametric foundations are briefly analyzed in the following:

1. As shown in Figure 1, the Fourier series multiscale solution has good overall convergence characteristics. Specifically, the composite Fourier series of the deflection of mid-surface, average rotations and bending moments of the thick plate converge well within the solution domain; the composite Fourier series of the deflection of mid-surface and average rotations of the thick plate converge well on the boundary and at the corners of the solution domain; the composite Fourier series of the bending moments of the thick plate do not converge so well on the boundary and at the corners of the solution domain.

2. As shown in Figures 2-5, the influencing factors have different effects on the convergence characteristics of the Fourier series multiscale solution. Specifically, the adjustment of computational parameters has obvious effect on the convergence of the Fourier series multiscale solution; the adjustment of thickness-to-length ratio and length-to-width ratio have unobvious effect on the convergence of the Fourier series multiscale solution; the adjustment of boundary conditions from CCCC boundary condition to SSSS boundary condition and FFFF boundary condition leads to obviously worse convergence for the composite Fourier series of the deflection of mid-surface, while it is uncertain for the convergence of the composite Fourier series of average rotations and bending moments of the thick plate.

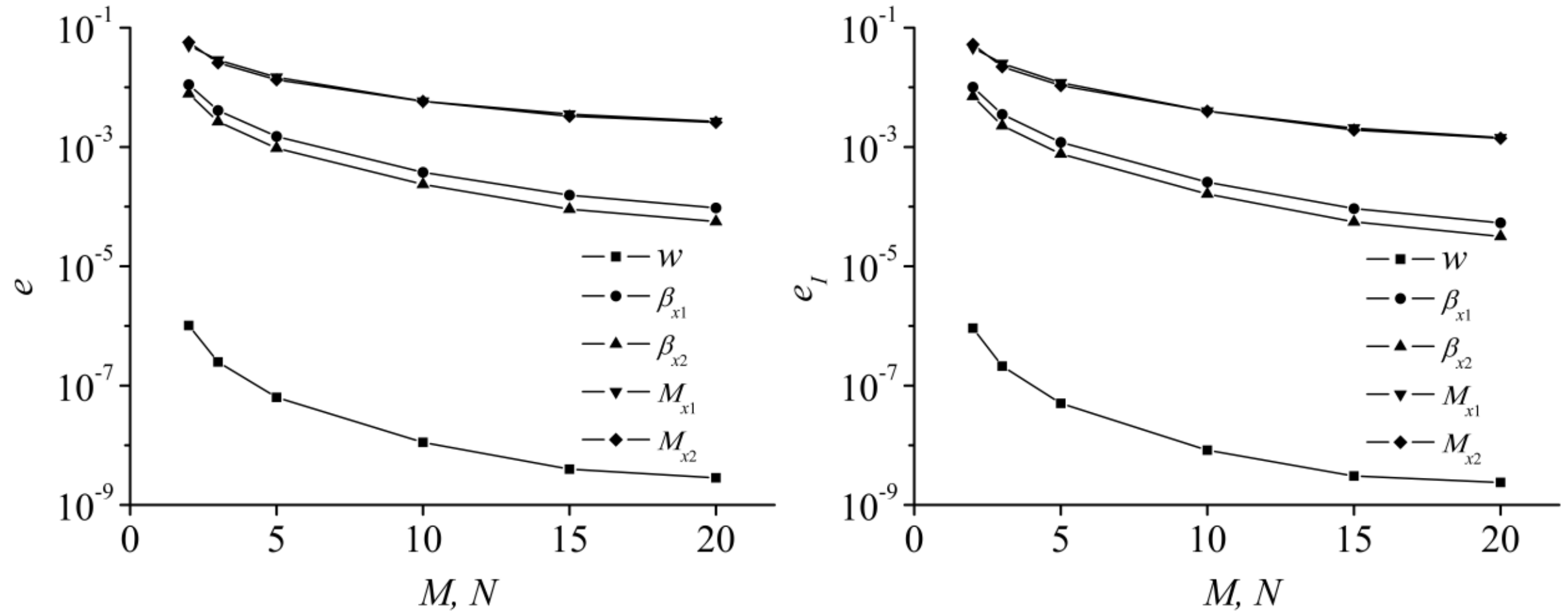

(a) (b)

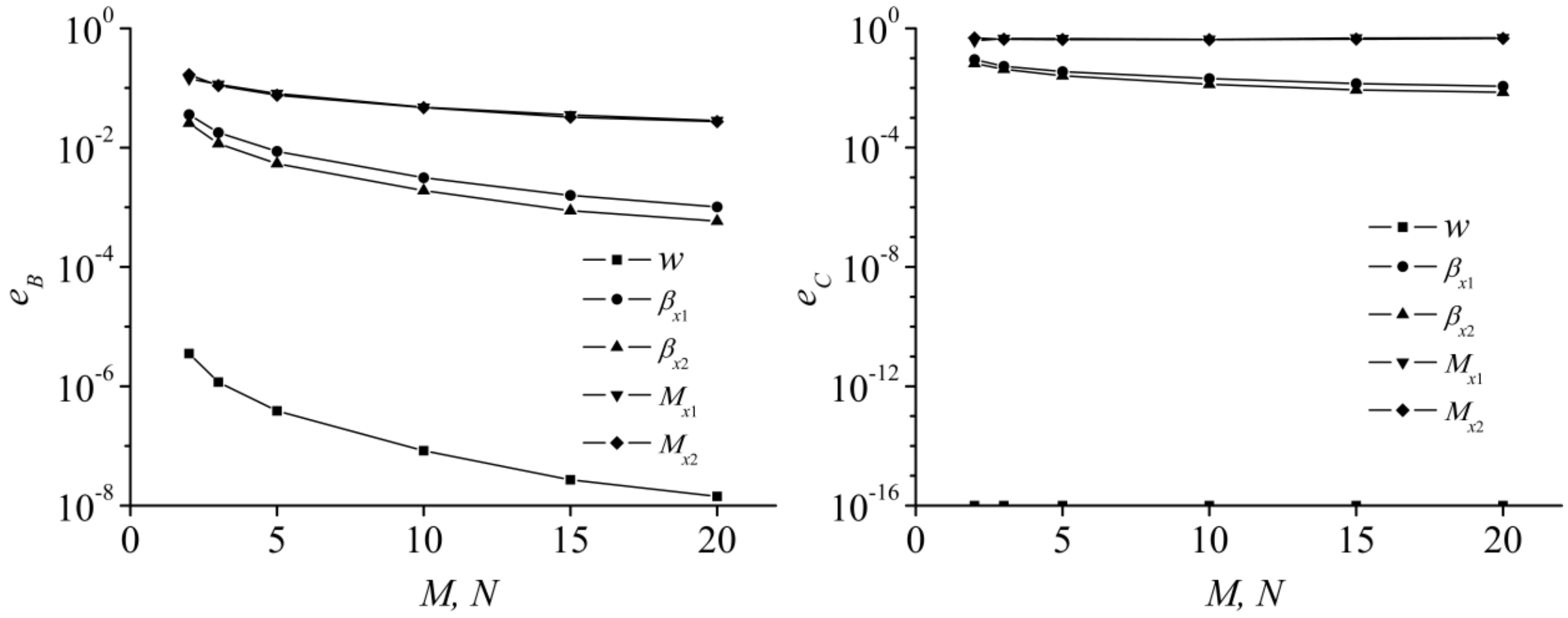

(c) (d)

Figure 1: Convergence characteristics of the Fourier series multiscale solution for the elastic bending of a thick plate resting on the biparametric foundation:

(a) $e$ - $M$ , $N$ curves, (b) $e_I$ - $M$ , $N$ curves, (c) $e_B$ - $M$ , $N$ curves, (d) $e_C$ - $M$ , $N$ curves.

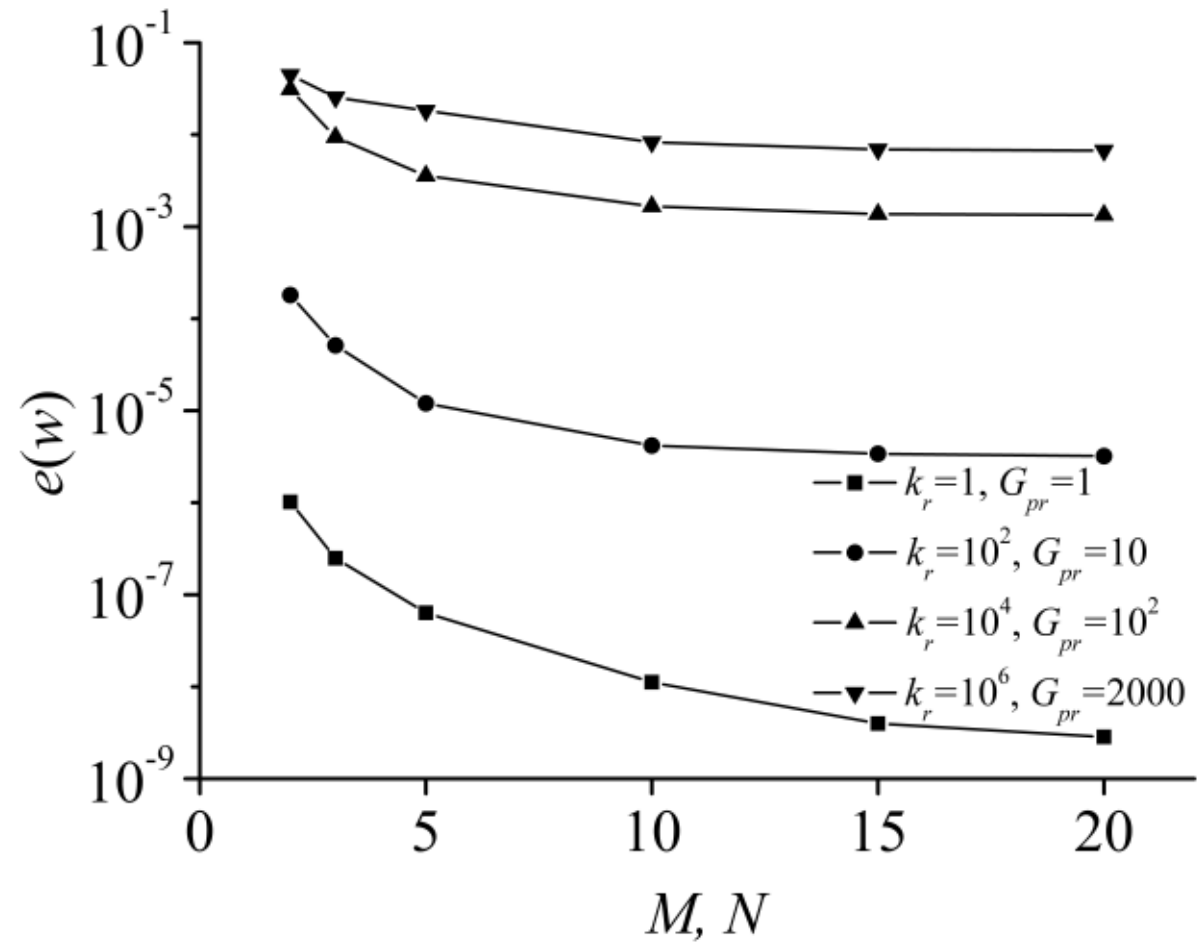

(a)

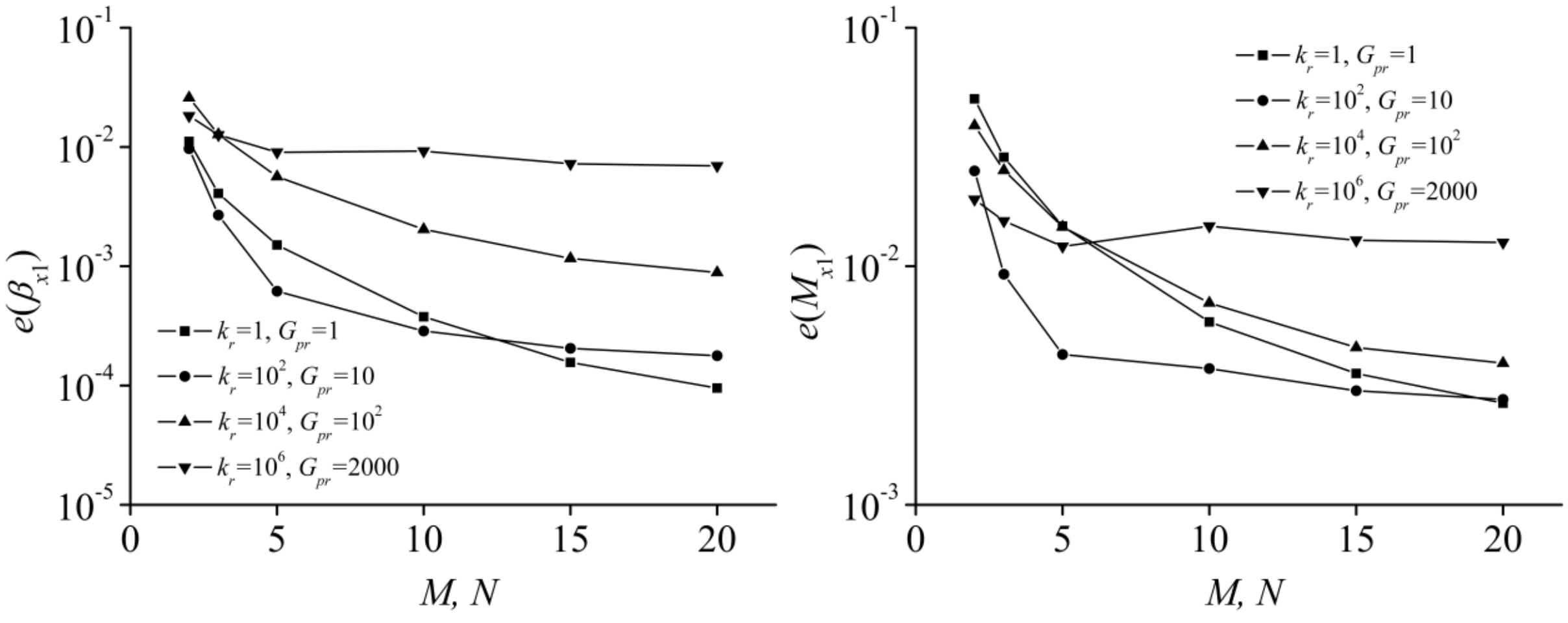

(b) (c)

Figure 2: Convergence comparison of the Fourier series multiscale solutions with different computational parameters $k_r$ and $G_{pr}$ :

(a) $e(w)$ - $M$ , $N$ curves, (b) $e(\beta_x)$ - $M$ , $N$ curves, (c) $e(M_x)$ - $M$ , $N$ curves.

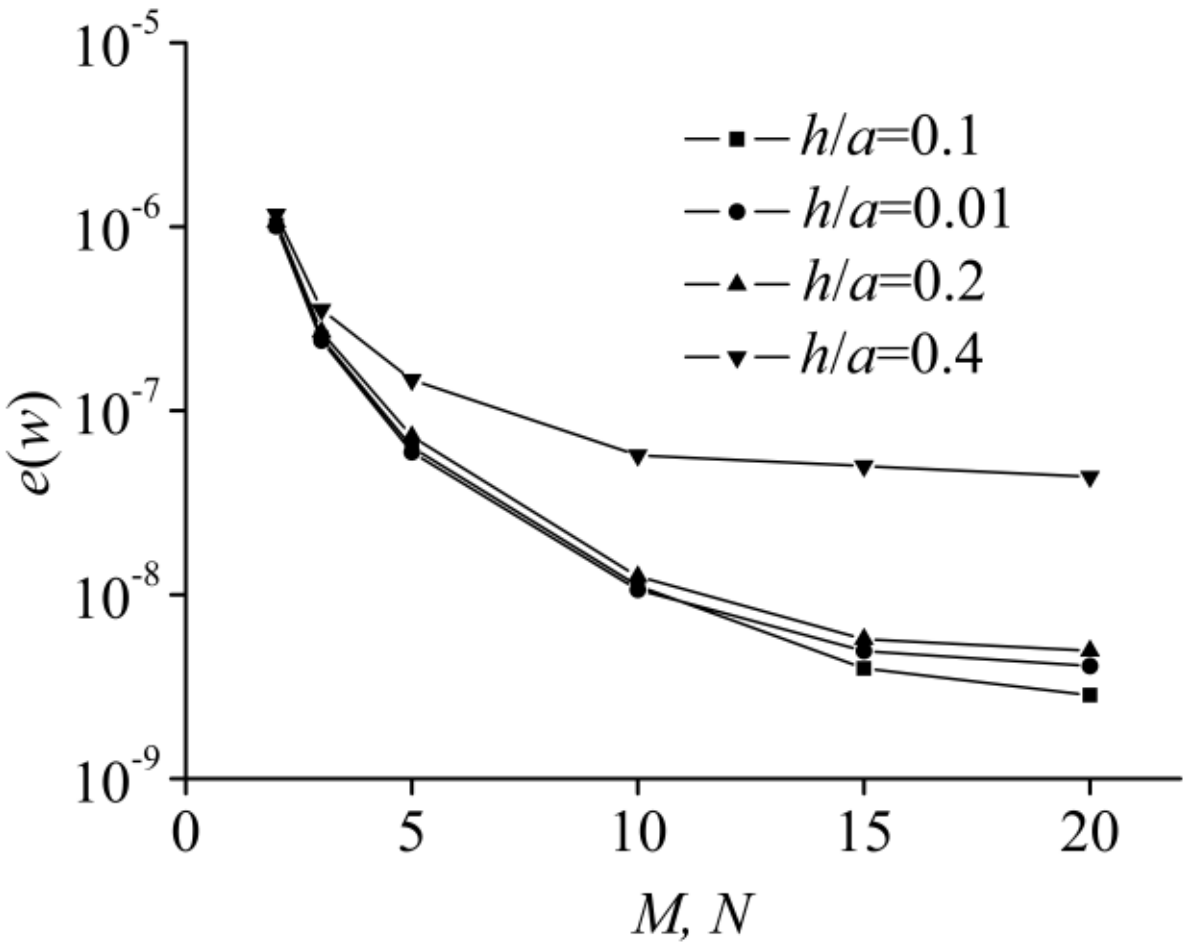

(a)

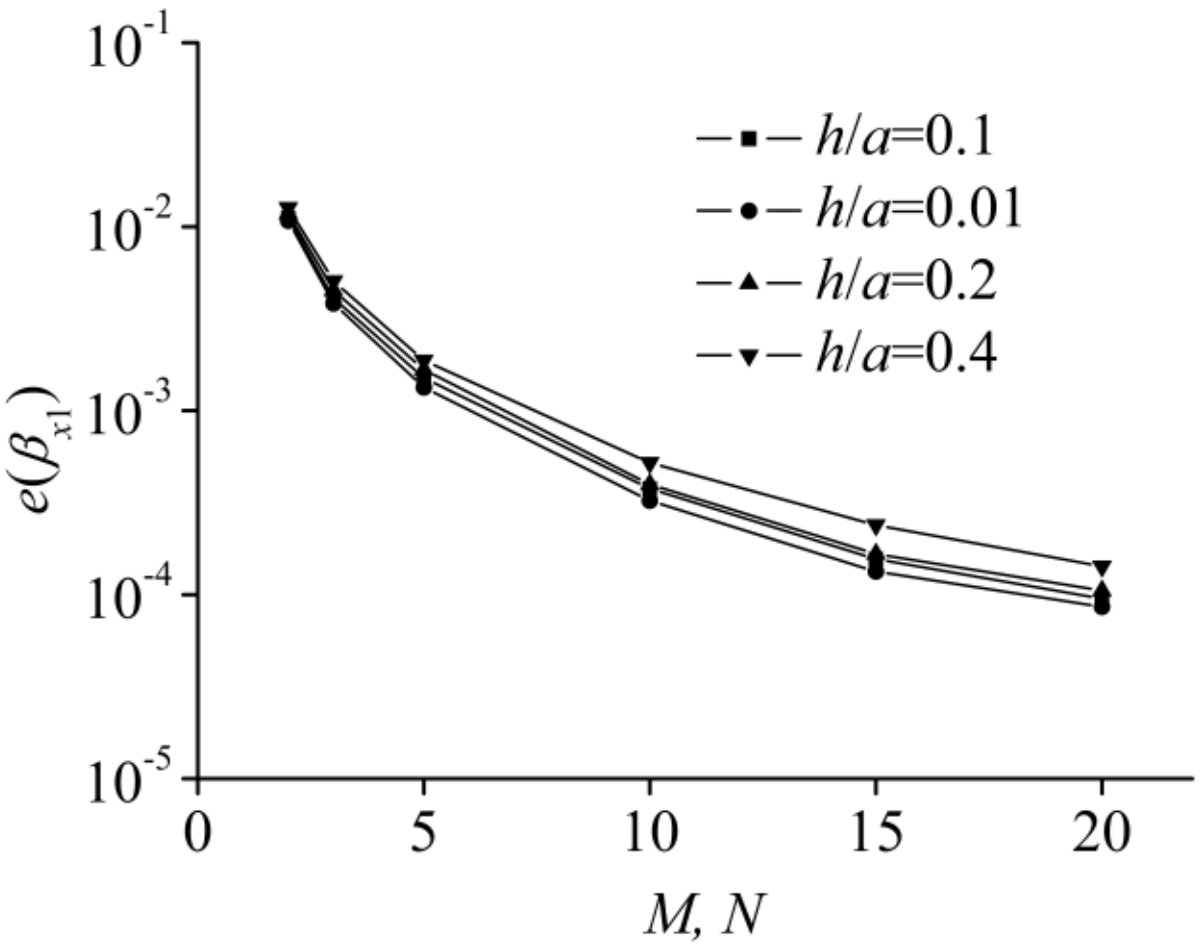

(b)

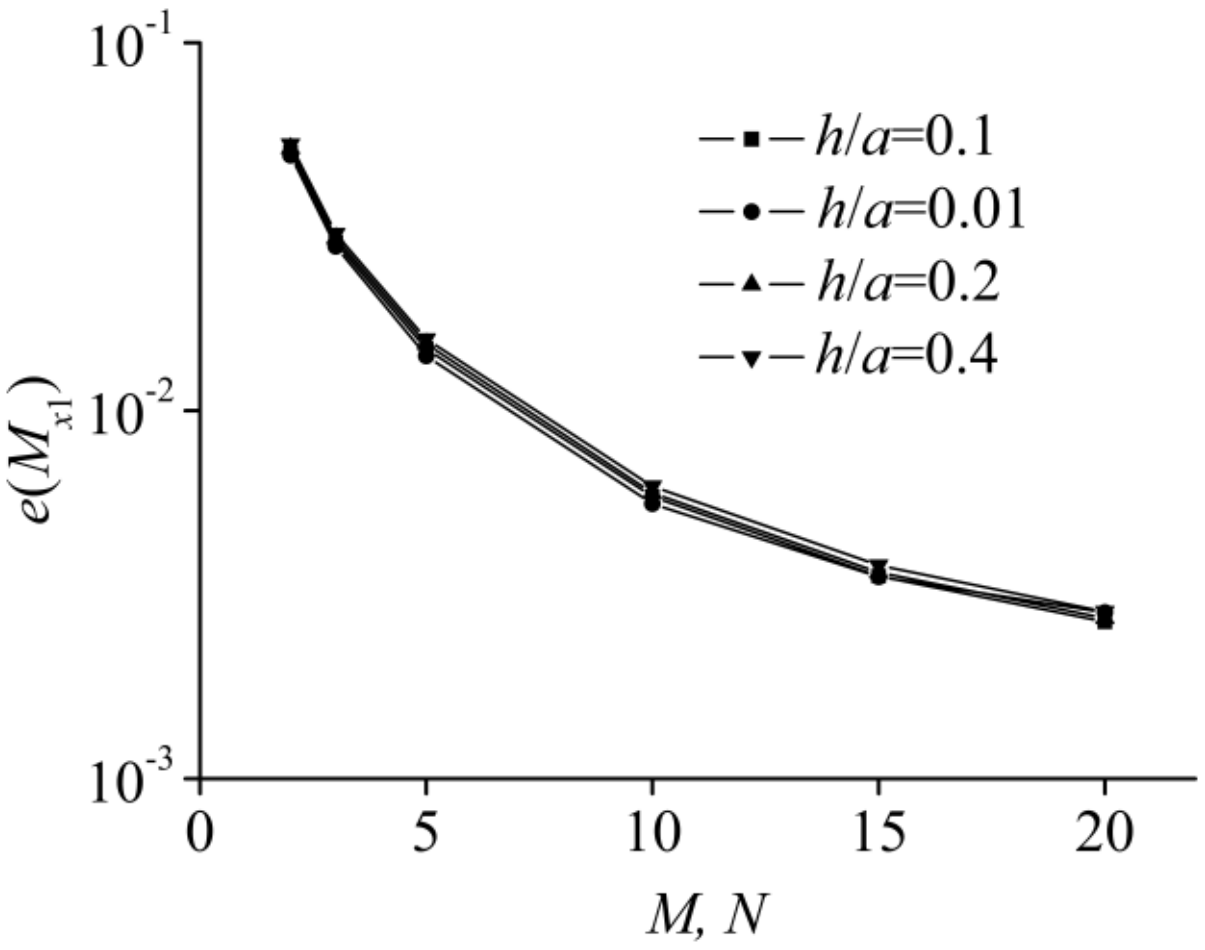

(c)

Figure 3: Convergence comparison of the Fourier series multiscale solutions with different thickness-to-length ratios:
(a) $e(w)$ - $M$ , $N$ curves, (b) $e(\beta_x)$ - $M$ , $N$ curves, (c) $e(M_x)$ - $M$ , $N$ curves.

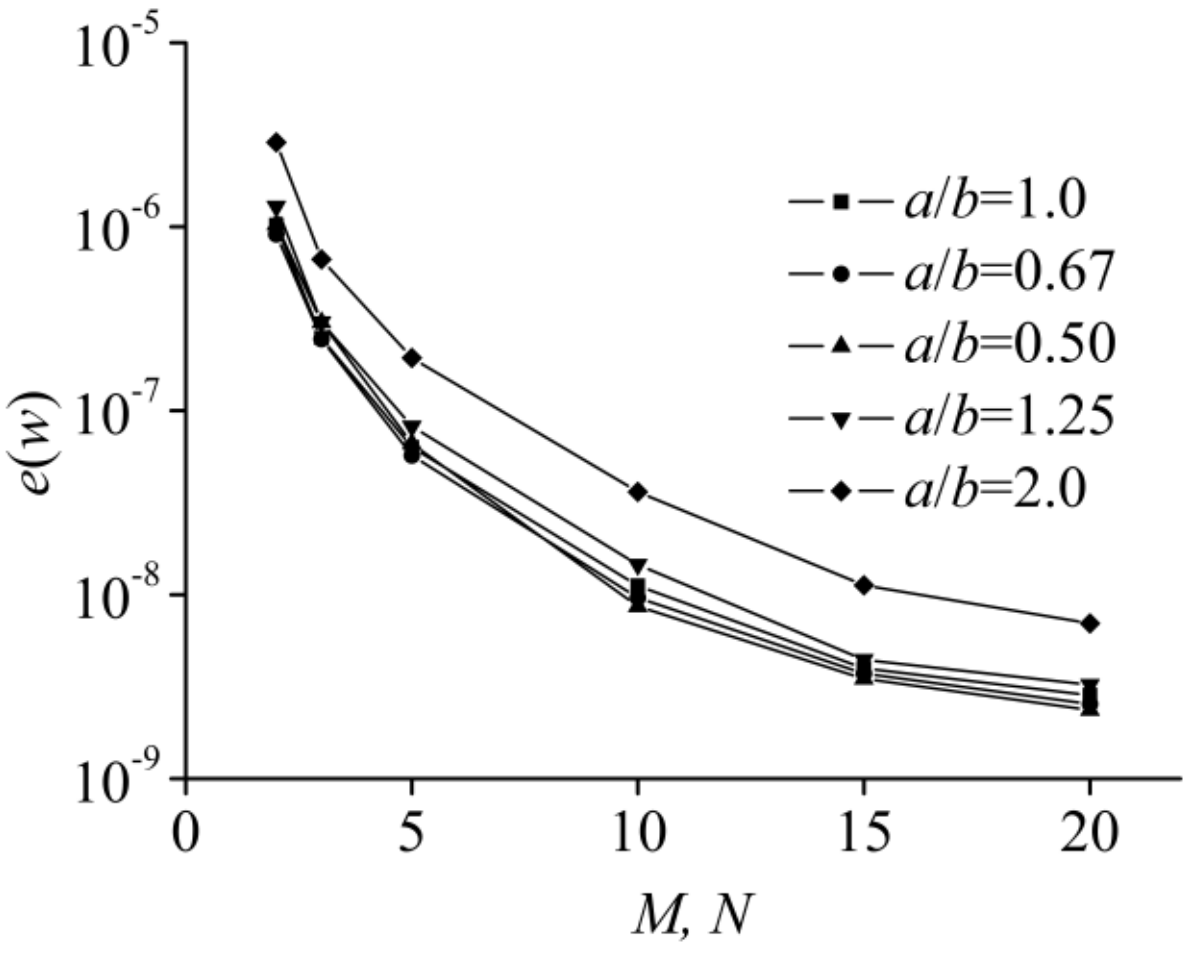

(a)

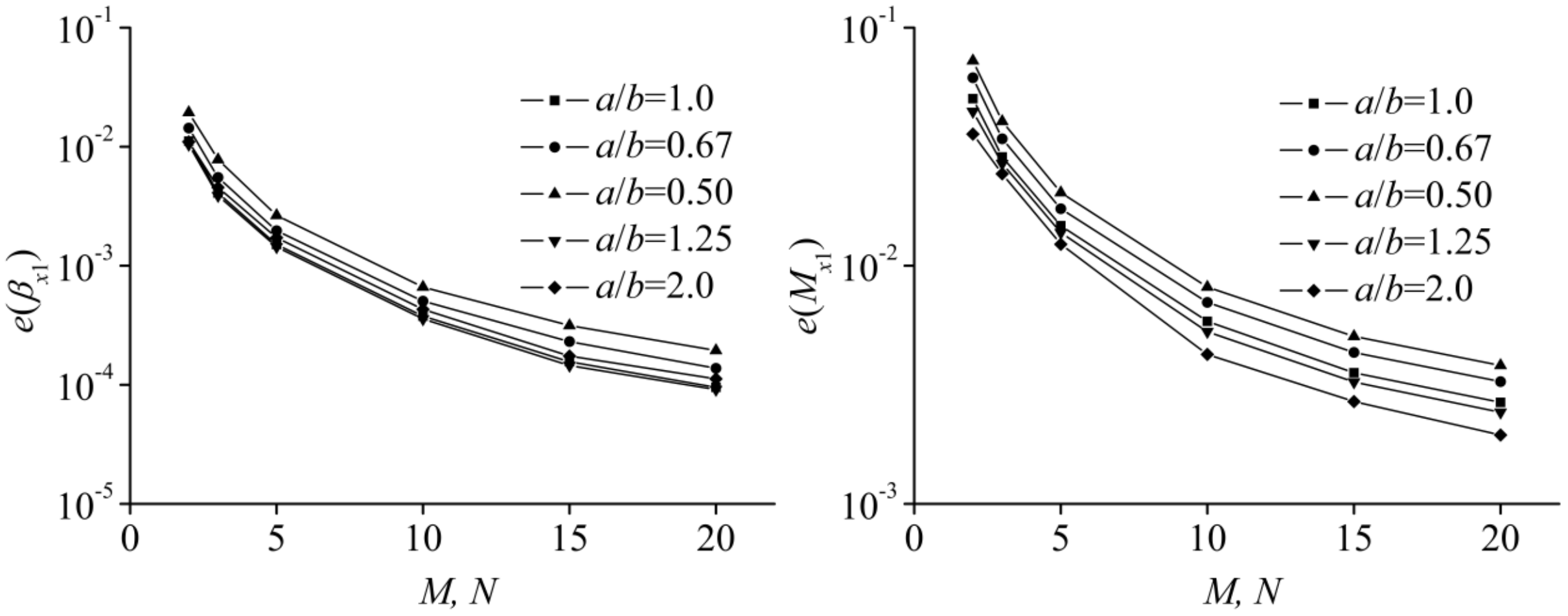

(b) (c)

Figure 4: Convergence comparison of the Fourier series multiscale solutions with different length-to-width ratios:
(a) $e(w)$-$M$, $N$ curves, (b) $e(\beta_x)$-$M$, $N$ curves, (c) $e(M_x)$-$M$, $N$ curves.

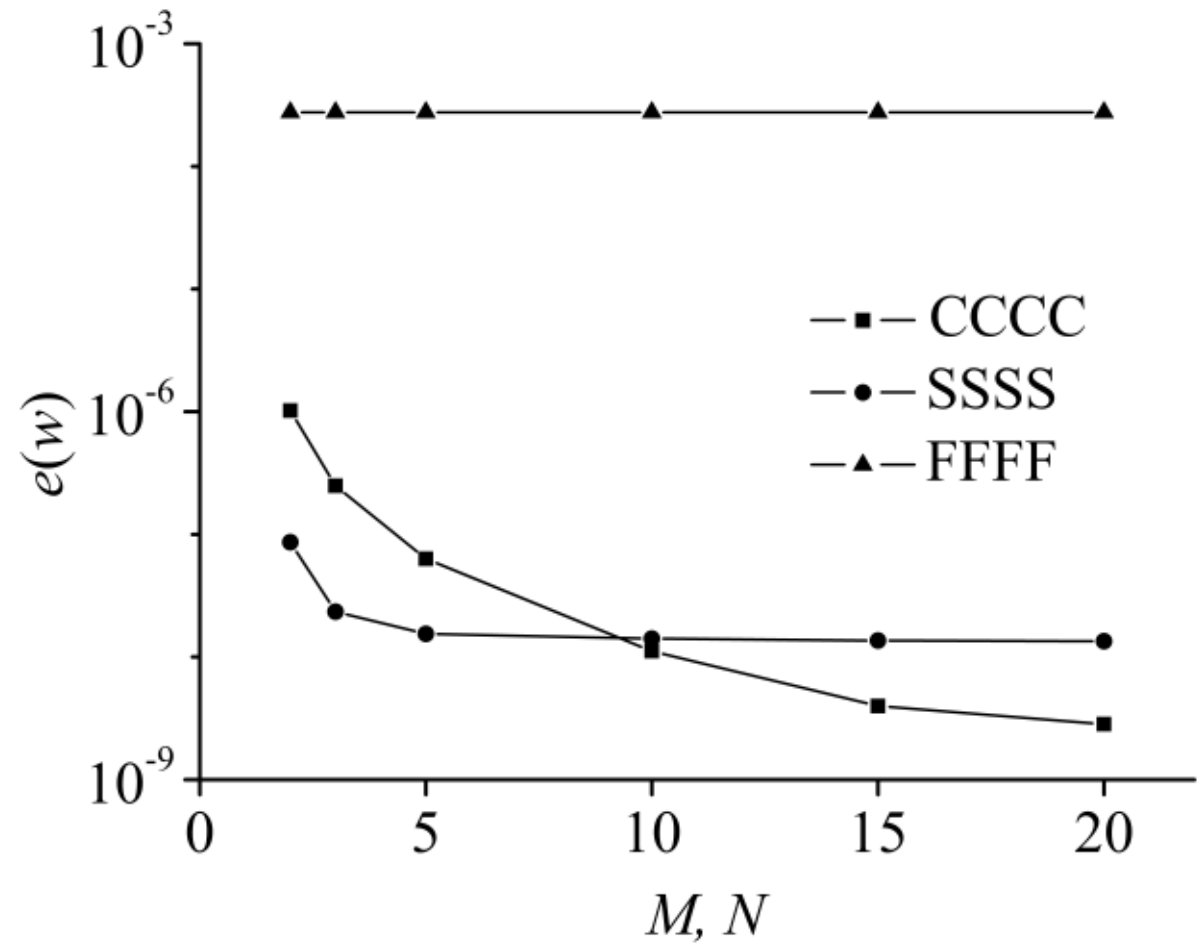

(a)

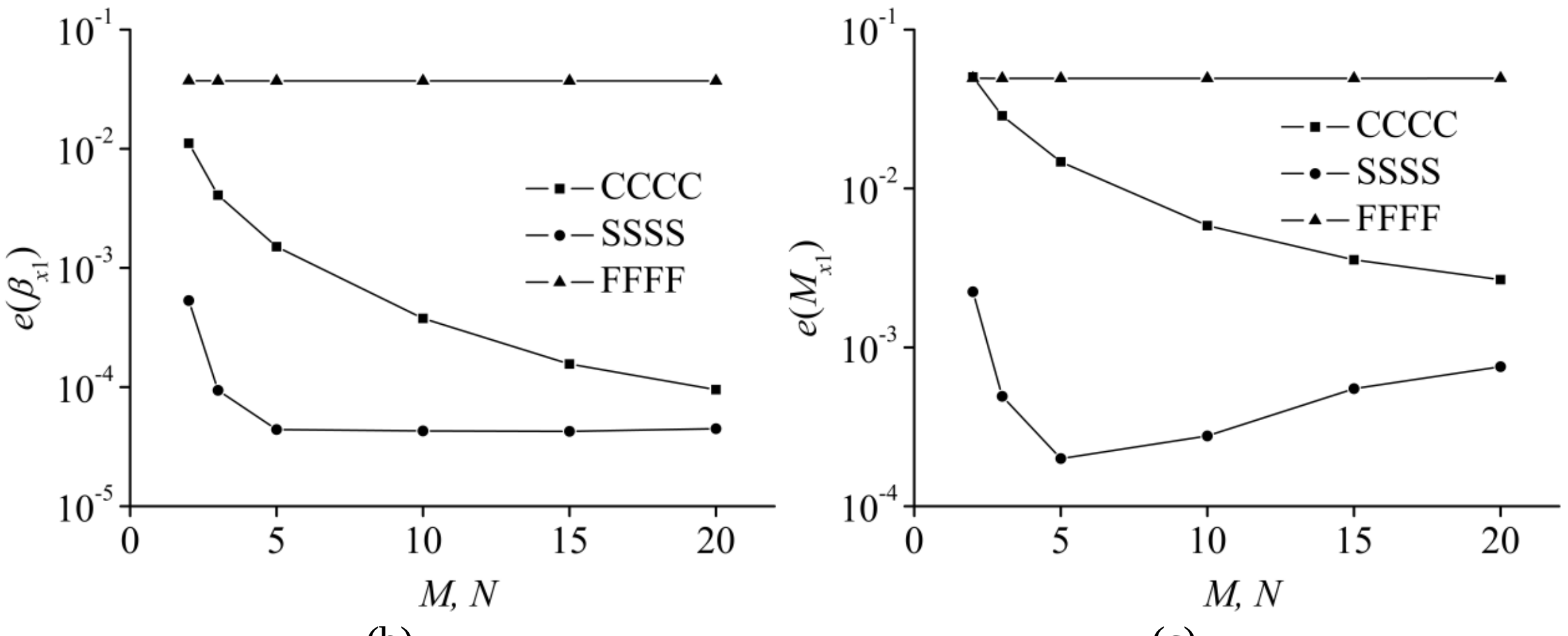

(b) (c)

Figure 5: Convergence comparison of the Fourier series multiscale solutions with different boundary conditions:
(a) $e(w)$-$M$, $N$ curves, (b) $e(\beta_x)$-$M$, $N$ curves, (c) $e(M_x)$-$M$, $N$ curves.

*5.2. Multiscale characteristics*

As shown in Table 2, based on the reference computational scheme in the first set of numerical experiment, we keep the computational parameter $k_r = 10^4$ for the elastic bending of thick plates on biparametric foundations unchanged and adjust the computational parameter $G_{pr}$ from 160 to 170, 180, 190 and 300 successively. We present in Figures 6-10 the corresponding reference results and the computed results of the Fourier series multiscale method with the number of truncated terms $M = N = 20$. It is shown that, during the adjustment process of computational parameters, the computed results with the Fourier series multiscale method exactly match the reference results, and the computed results within the solution domain are charactered by periodical fluctuation or a typical multiscale phenomenon, the boundary layer.

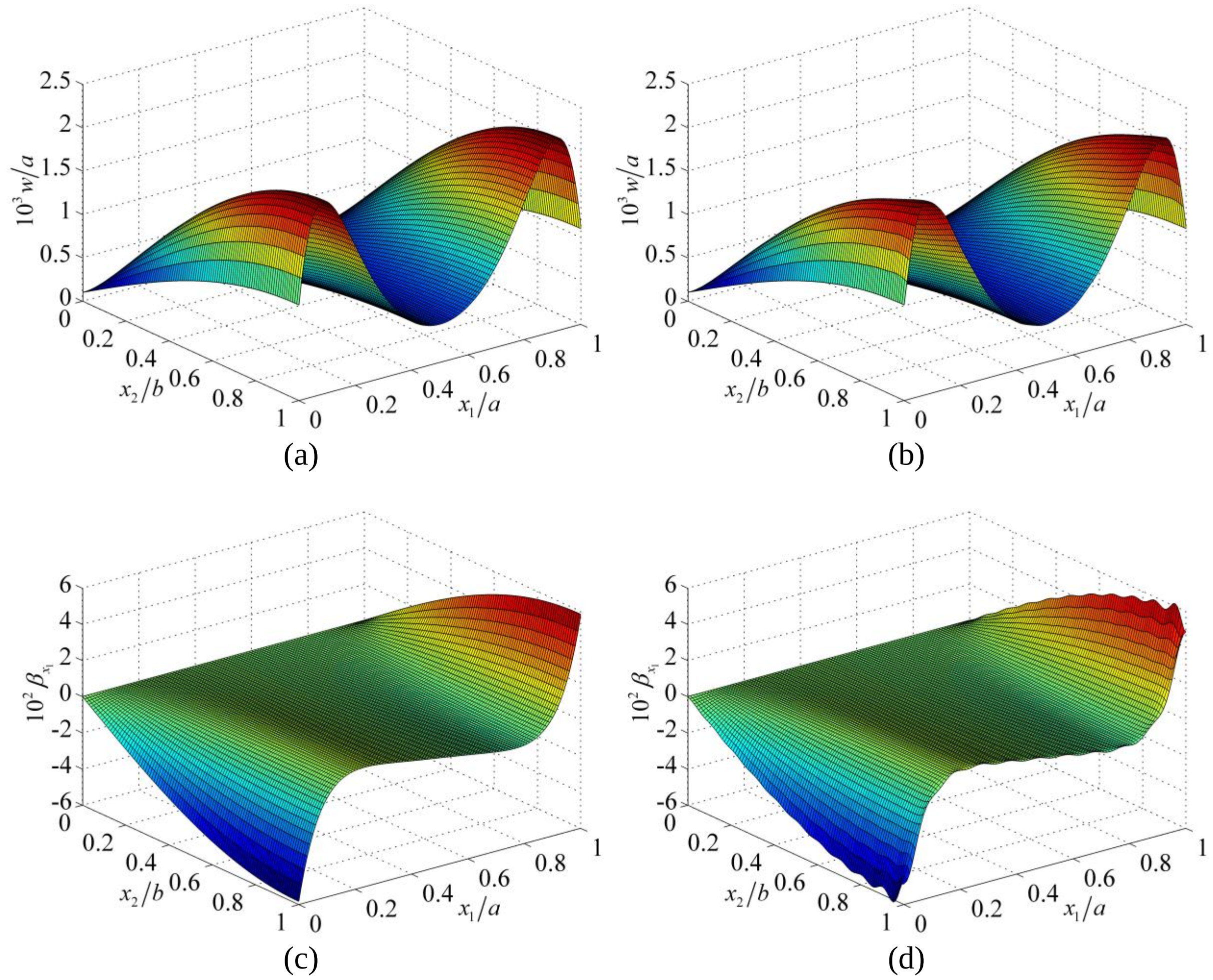

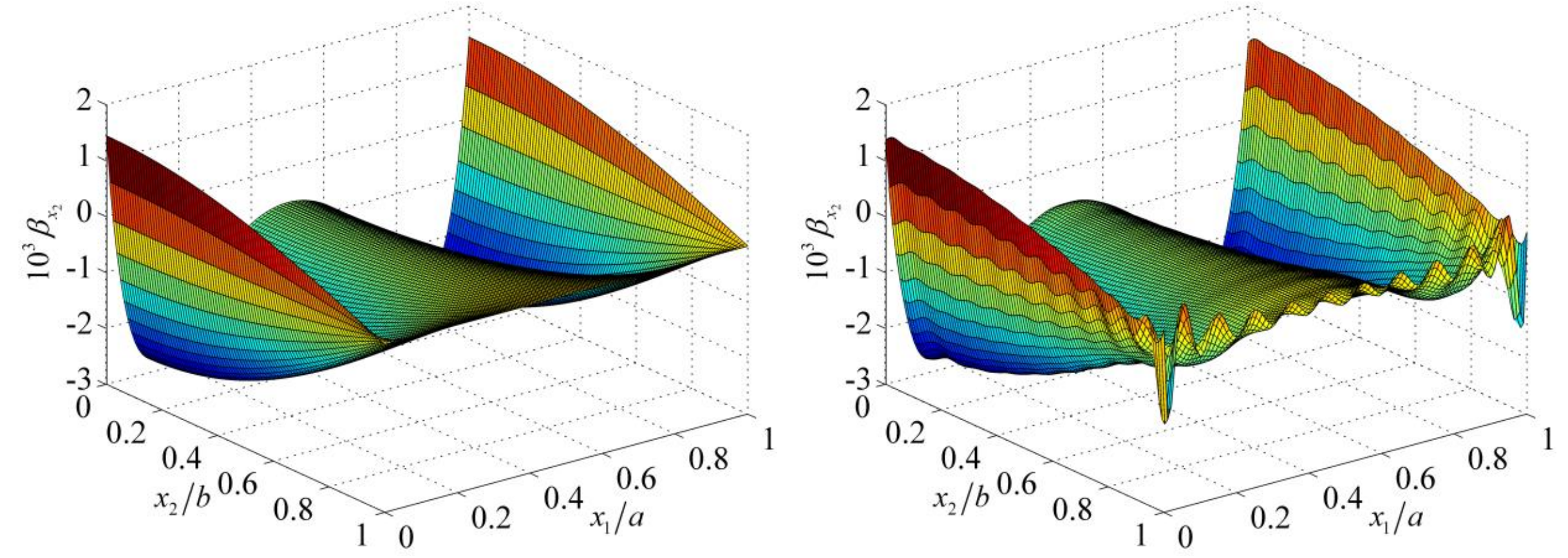

(e) (f)

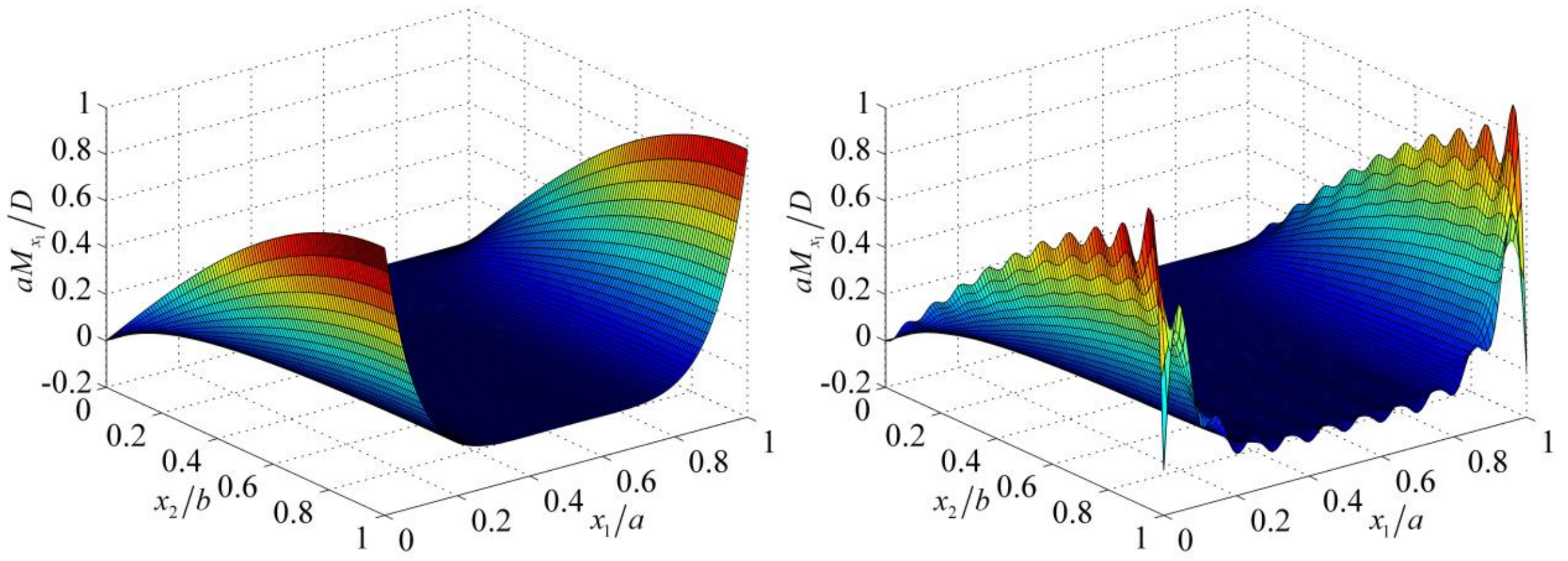

(g) (h)

Figure 6: The deflection, rotation and moment distributions for $k_r = 10^4$ and $G_{pr} = 160$:
(a) $w$ (reference), (b) $w$ (calculated), (c) $\beta_{x_1}$ (reference), (d) $\beta_{x_1}$ (calculated),
(e) $\beta_{x_2}$ (reference), (f) $\beta_{x_2}$ (calculated), (g) $M_{x_1}$ (reference), (h) $M_{x_1}$ (calculated).

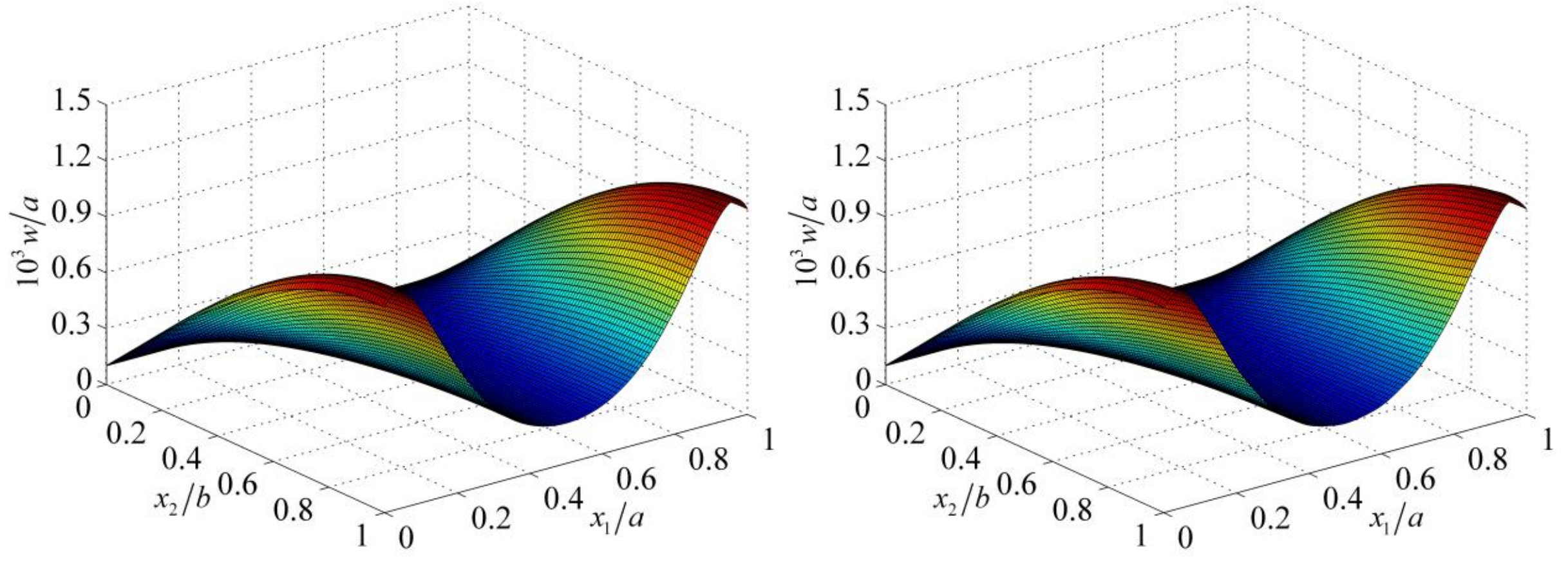

(a) (b)

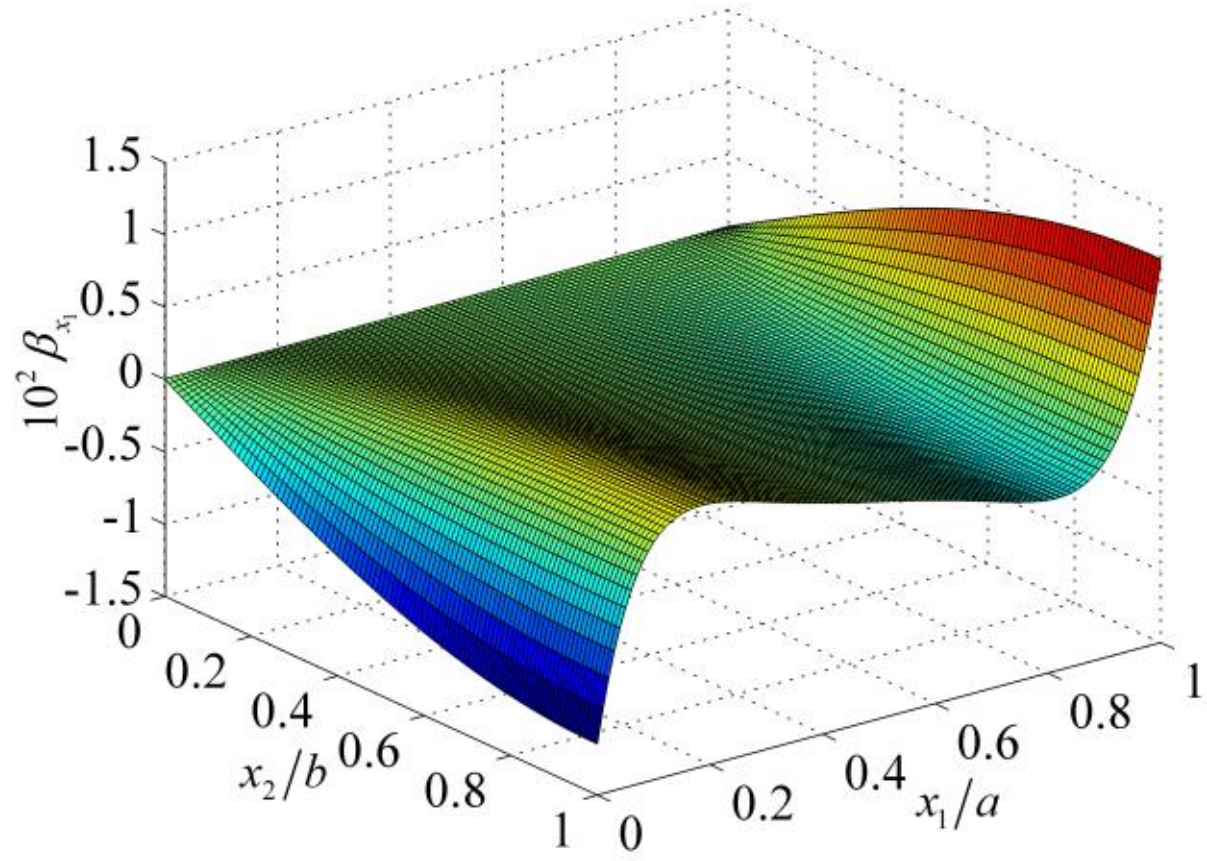

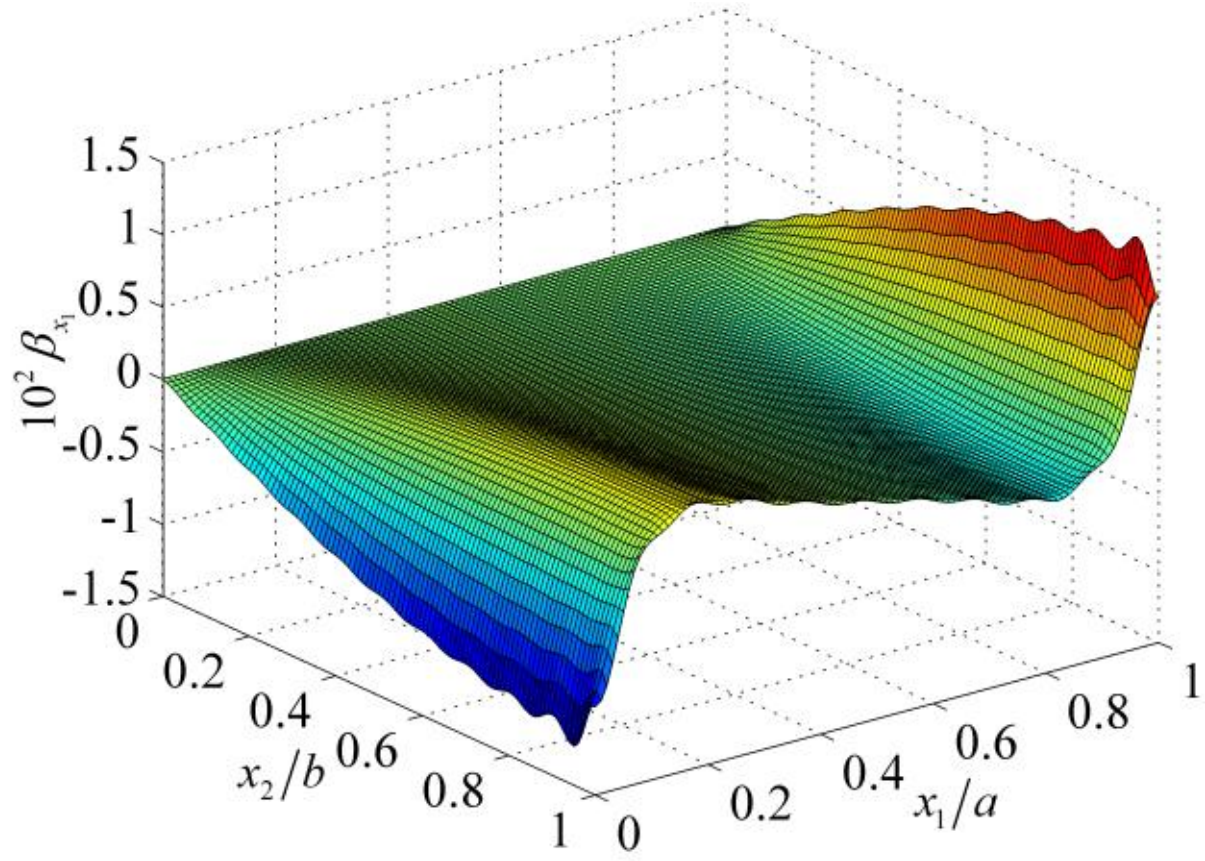

(c) (d)

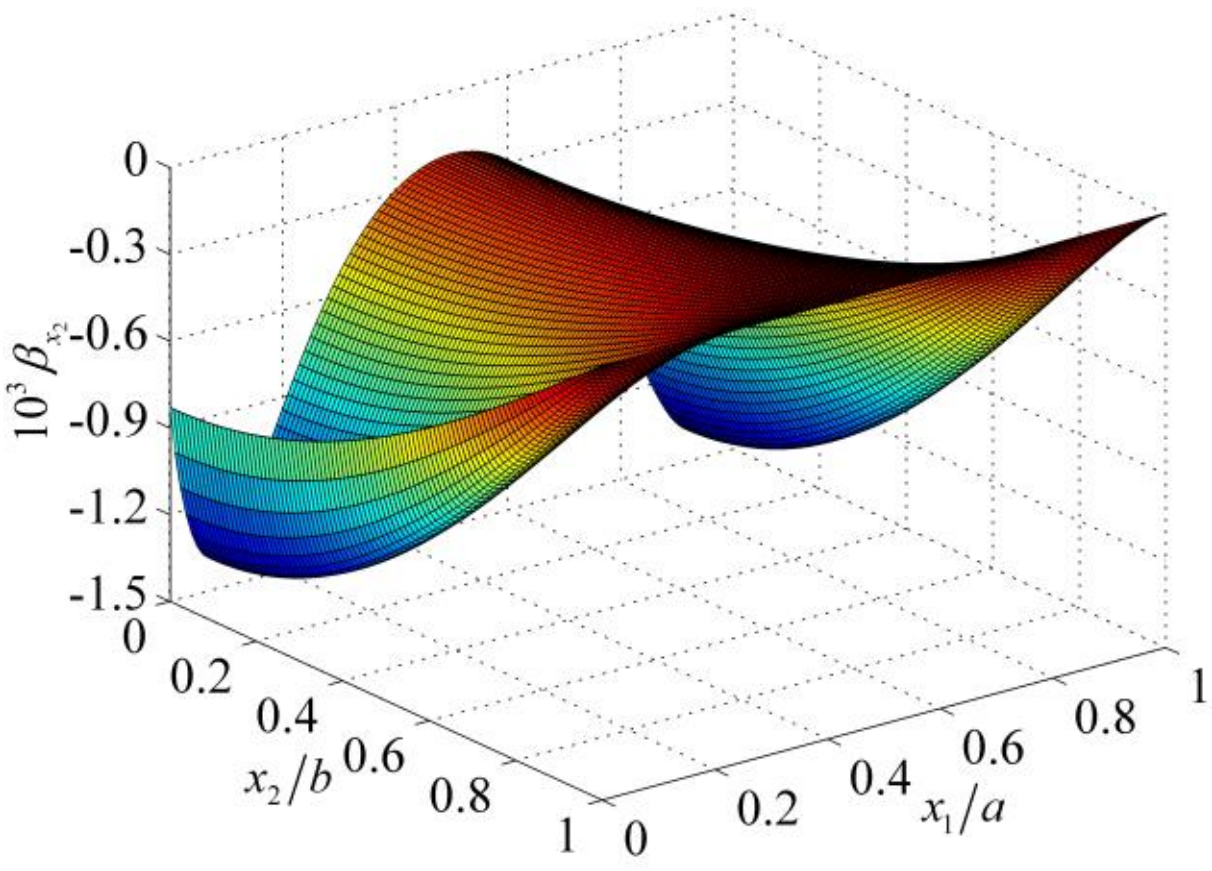

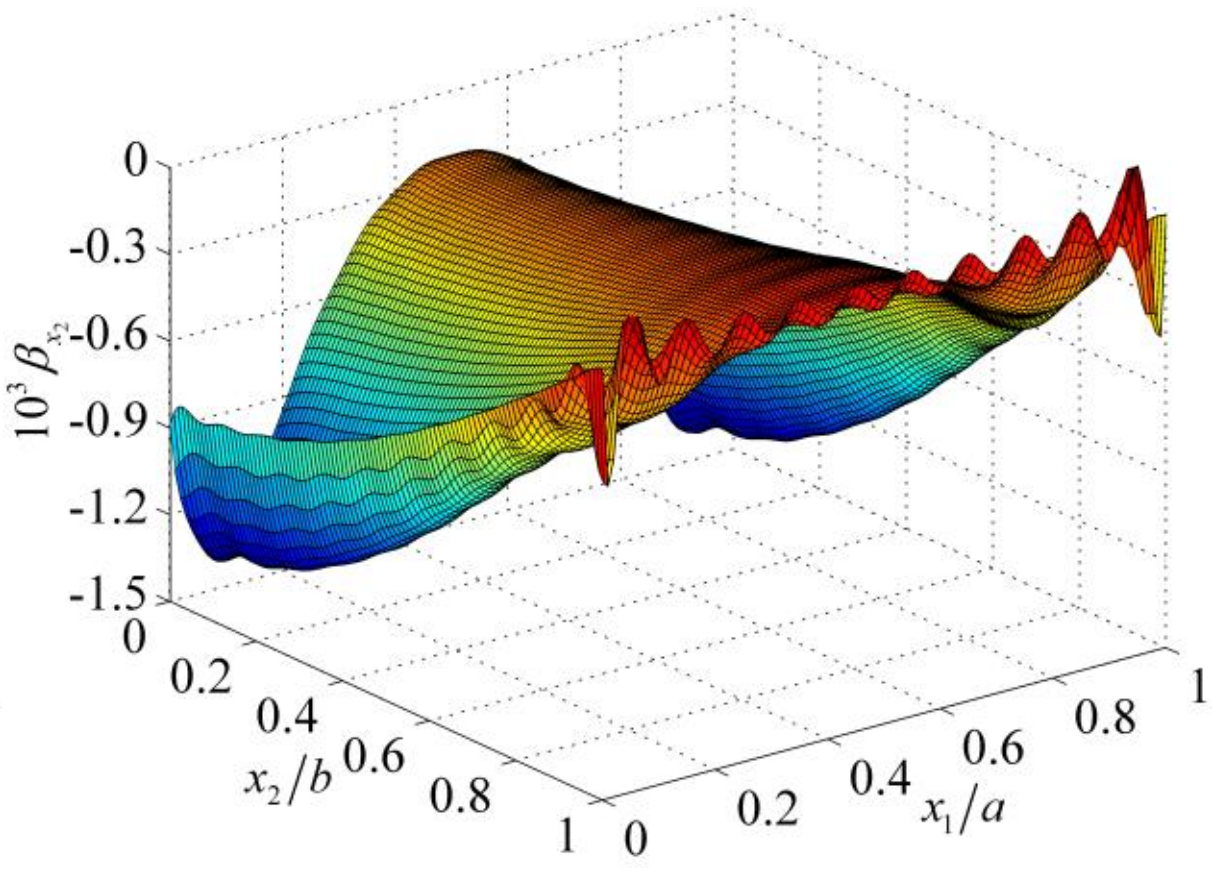

(e) (f)

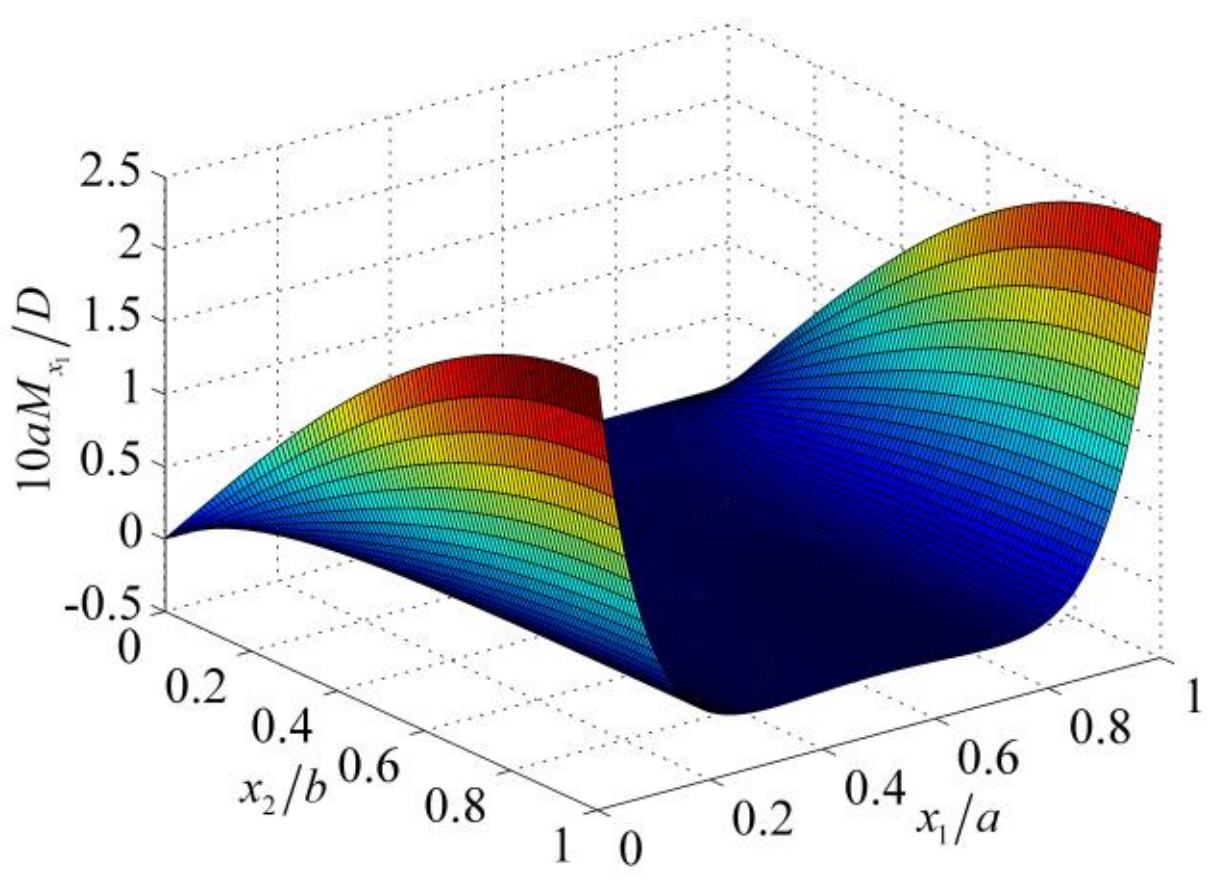

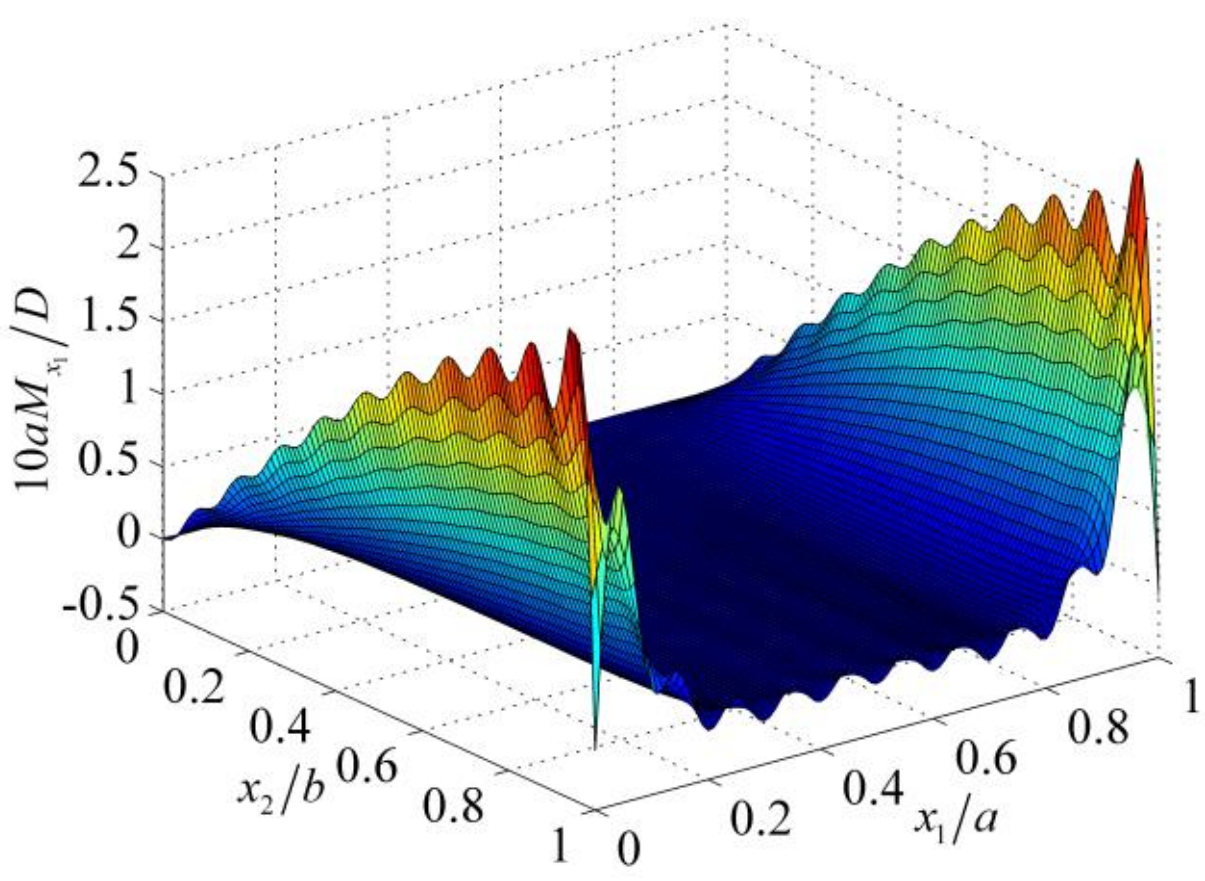

(g) (h)

Figure 7: The deflection, rotation and moment distributions for $k_r = 10^4$ and $G_{pr} = 170$: (a) $w$ (reference), (b) $w$ (calculated), (c) $\beta_{x_1}$ (reference), (d) $\beta_{x_1}$ (calculated), (e) $\beta_{x_2}$ (reference), (f) $\beta_{x_2}$ (calculated), (g) $M_{x_1}$ (reference), (h) $M_{x_1}$ (calculated).

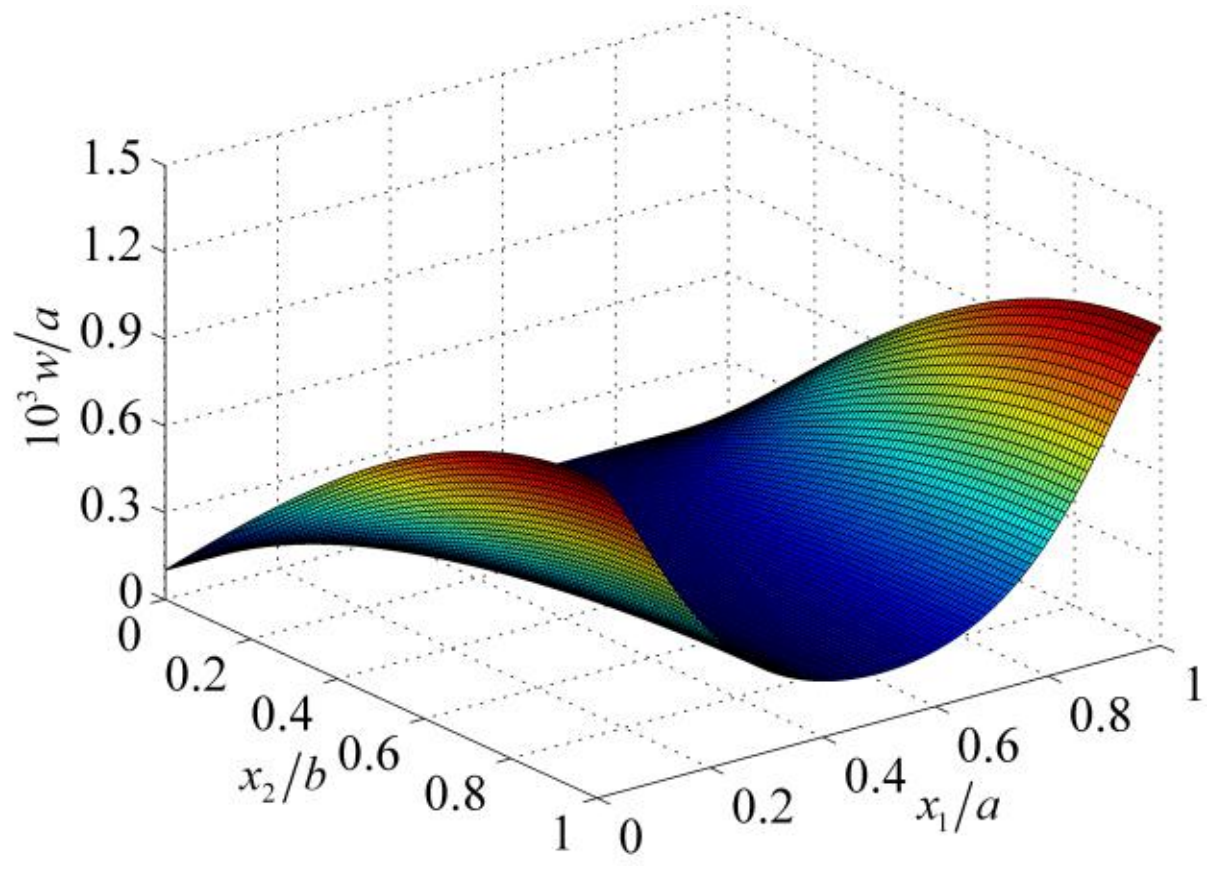

(a)

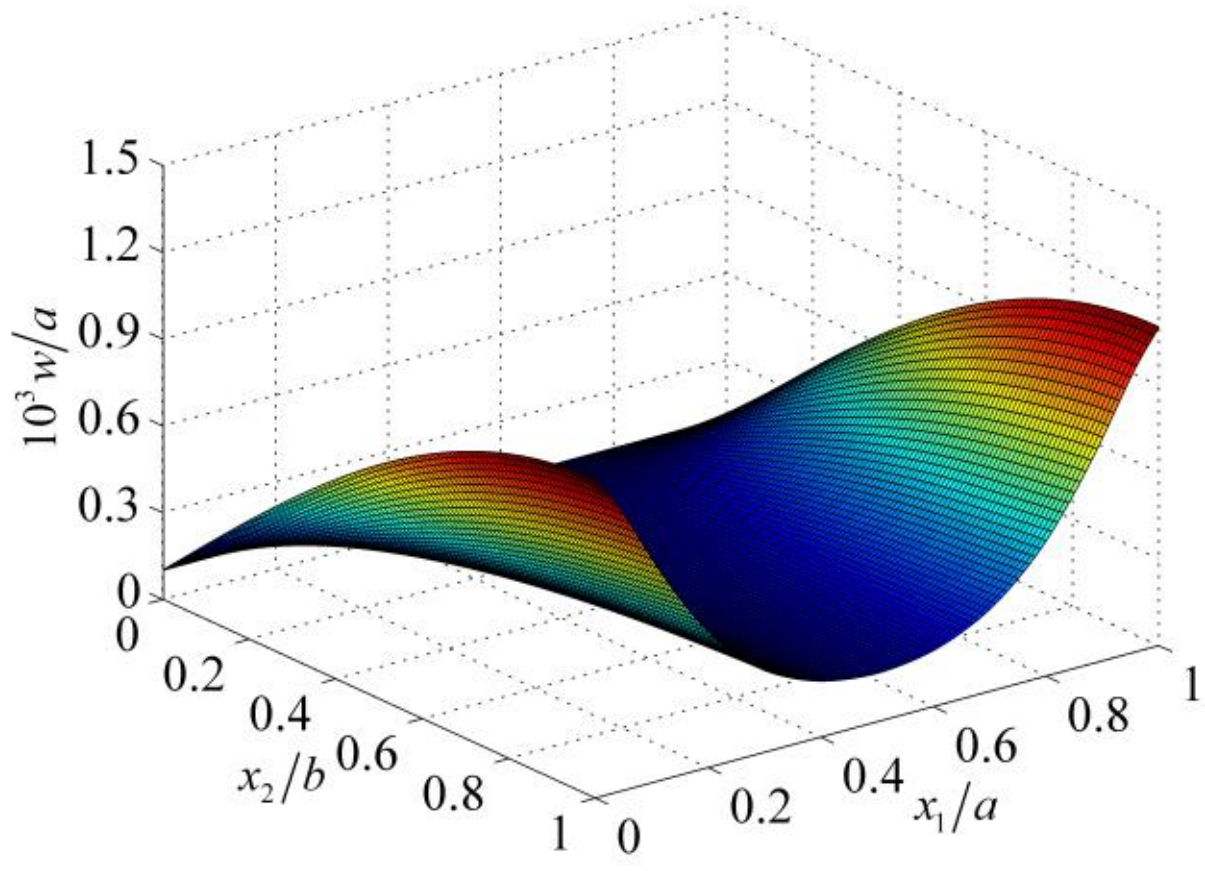

(b)

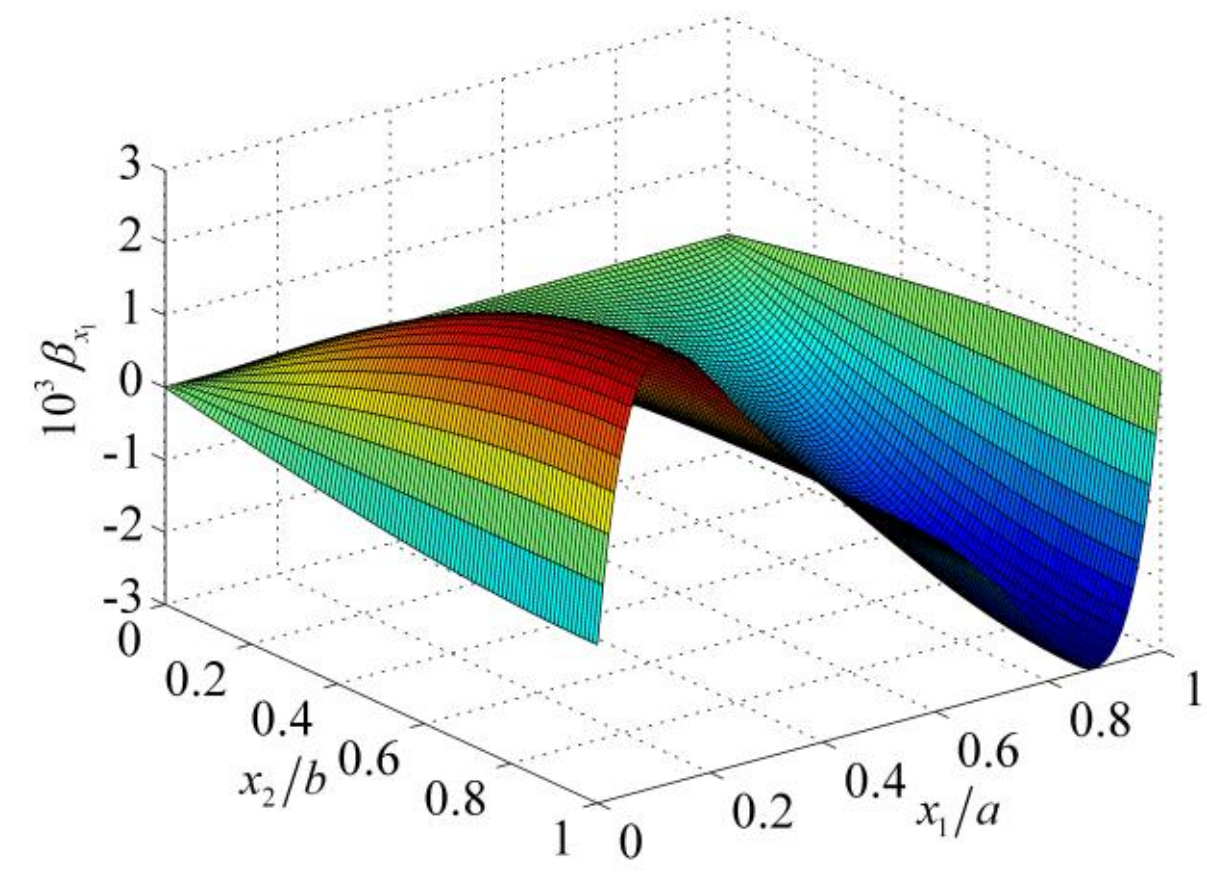

©

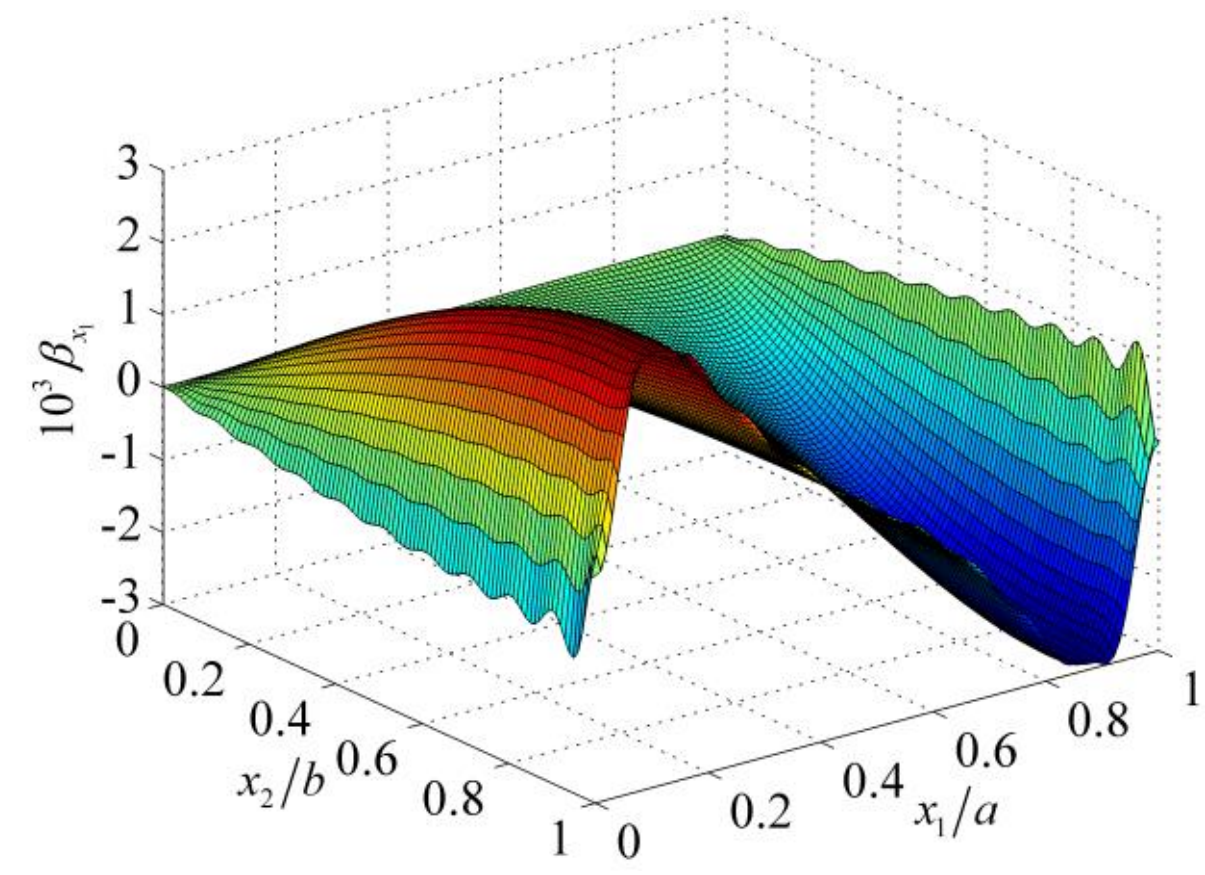

(d)

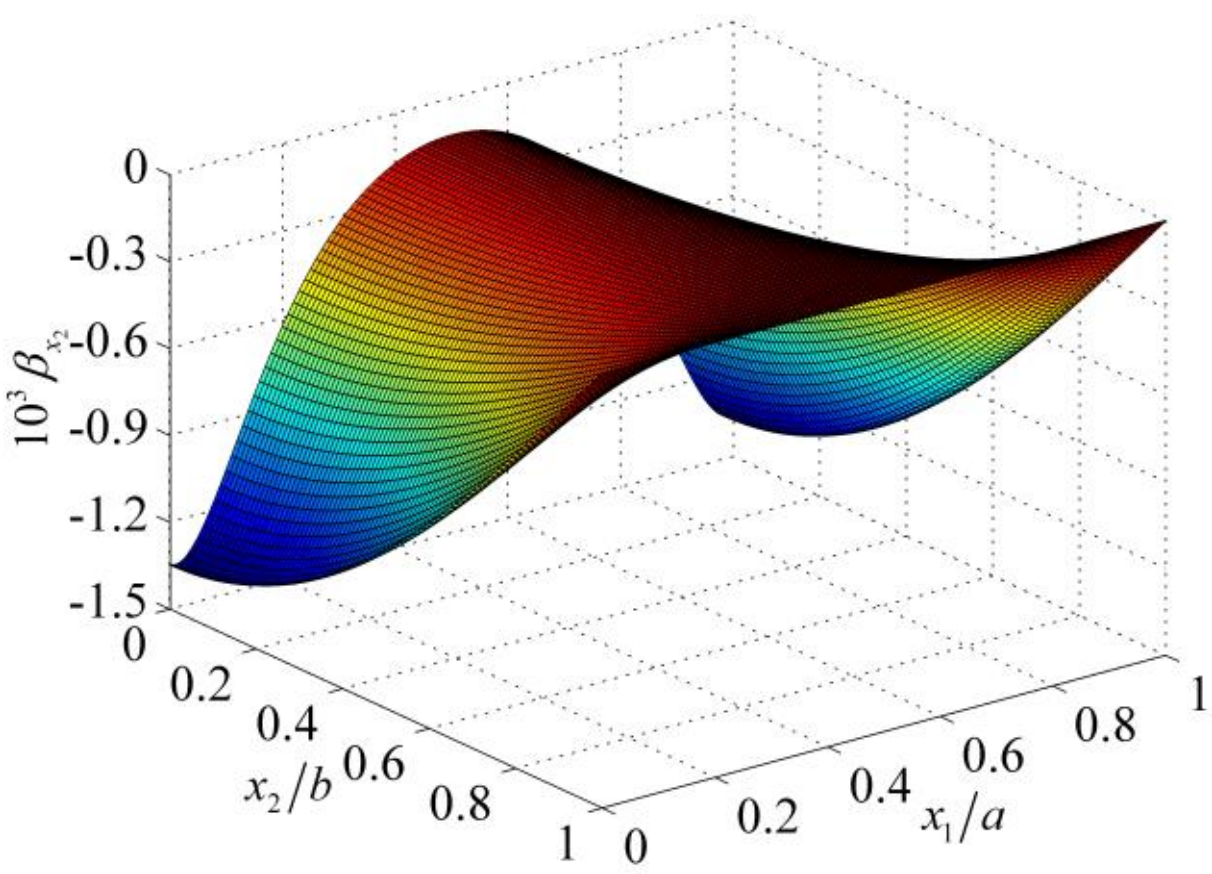

(e)

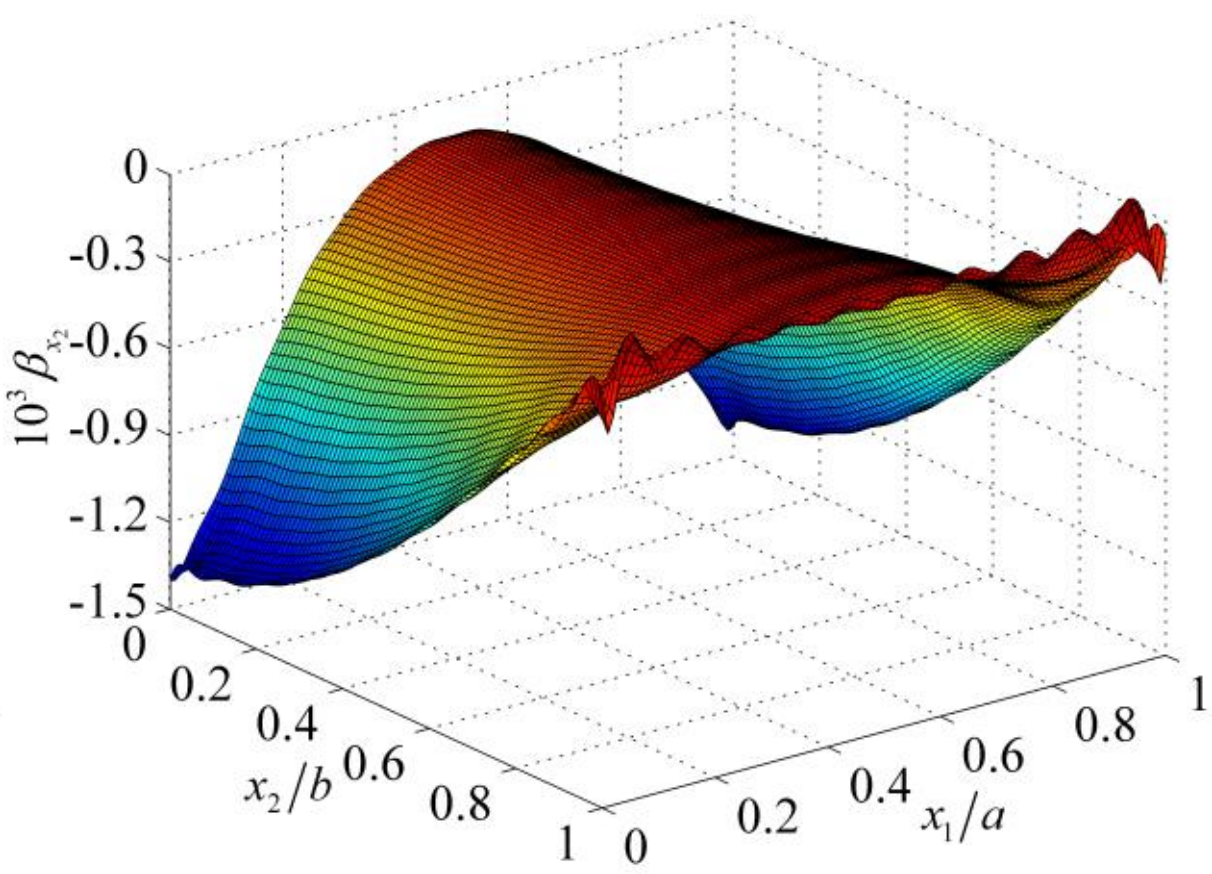

(f)

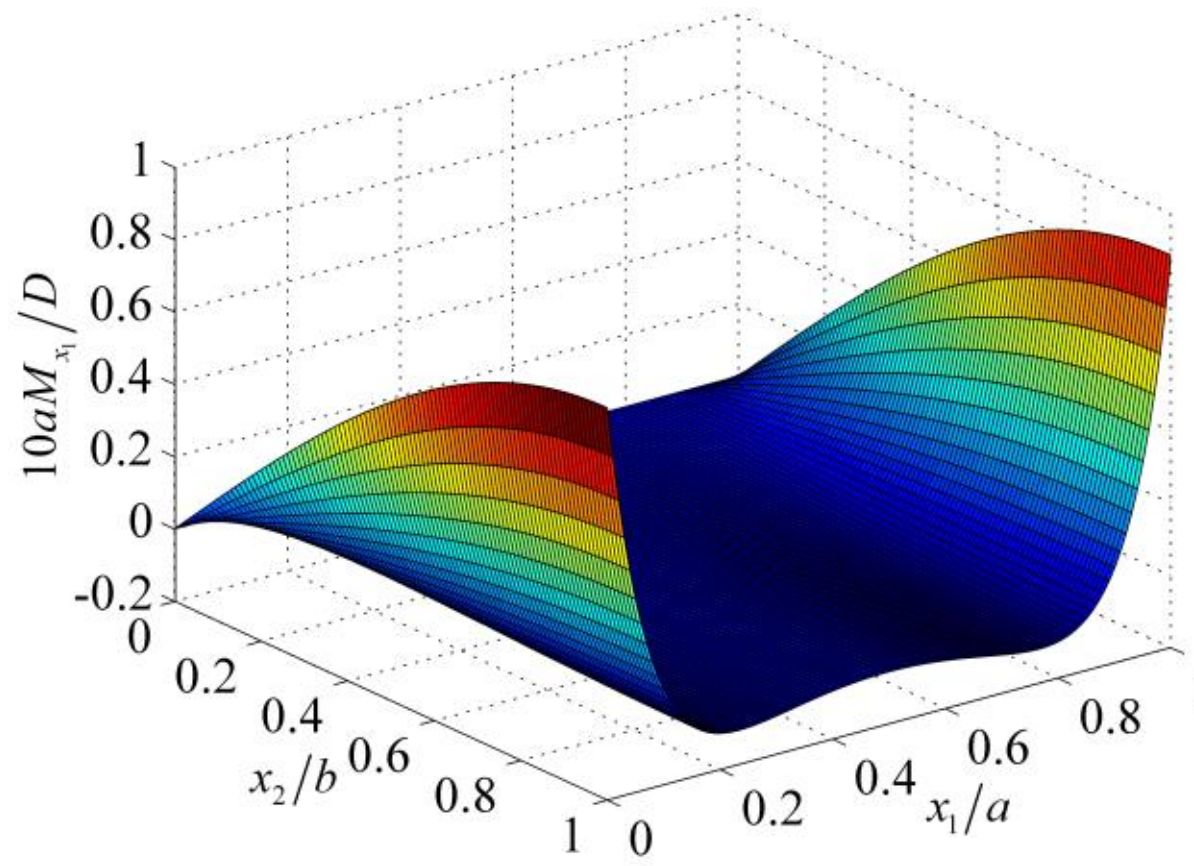

(g)

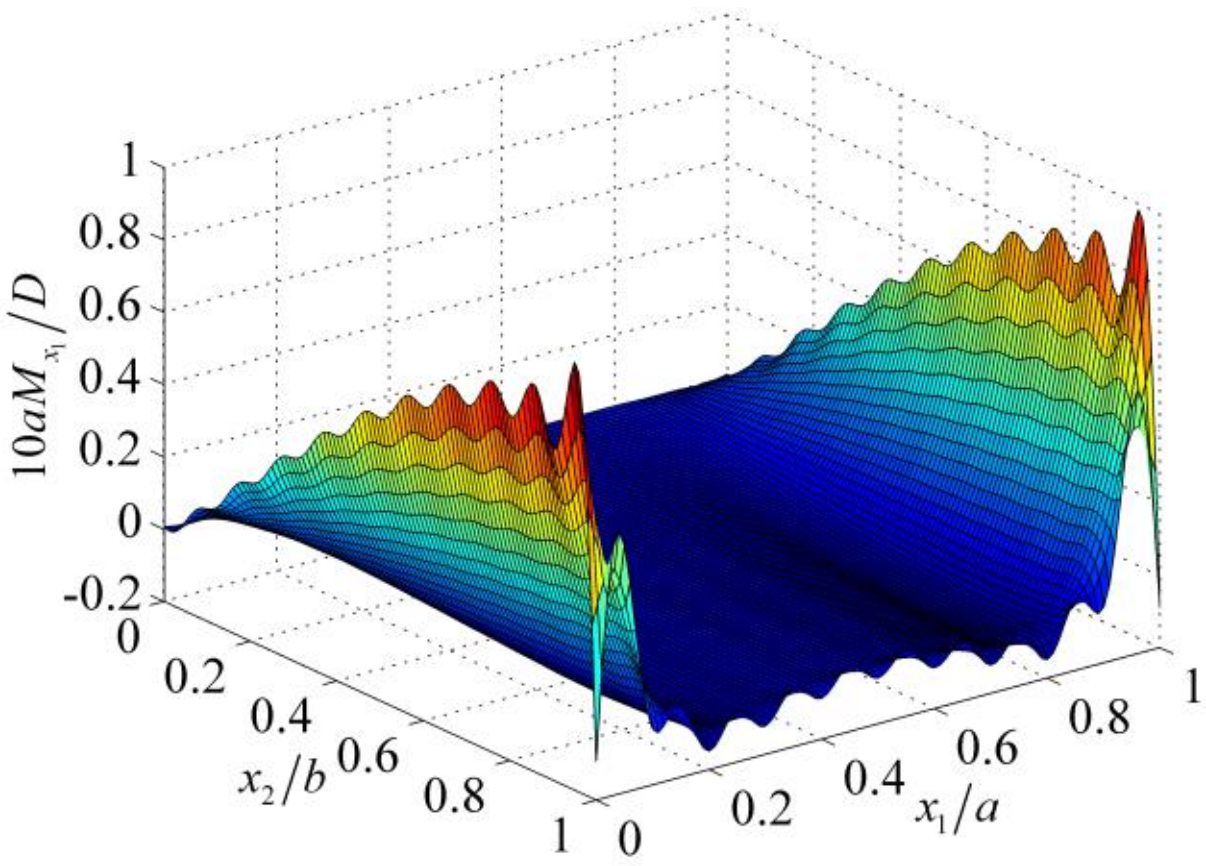

(h)

Figure 8: The deflection, rotation and moment distributions for $k_r = 10^4$ and $G_{pr} = 180$: (a) $w$ (reference), (b) $w$ (calculated), (c) $\beta_{x_1}$ (reference), (d) $\beta_{x_1}$ (calculated), (e) $\beta_{x_2}$ (reference), (f) $\beta_{x_2}$ (calculated), (g) $M_{x_1}$ (reference), (h) $M_{x_1}$ (calculated).

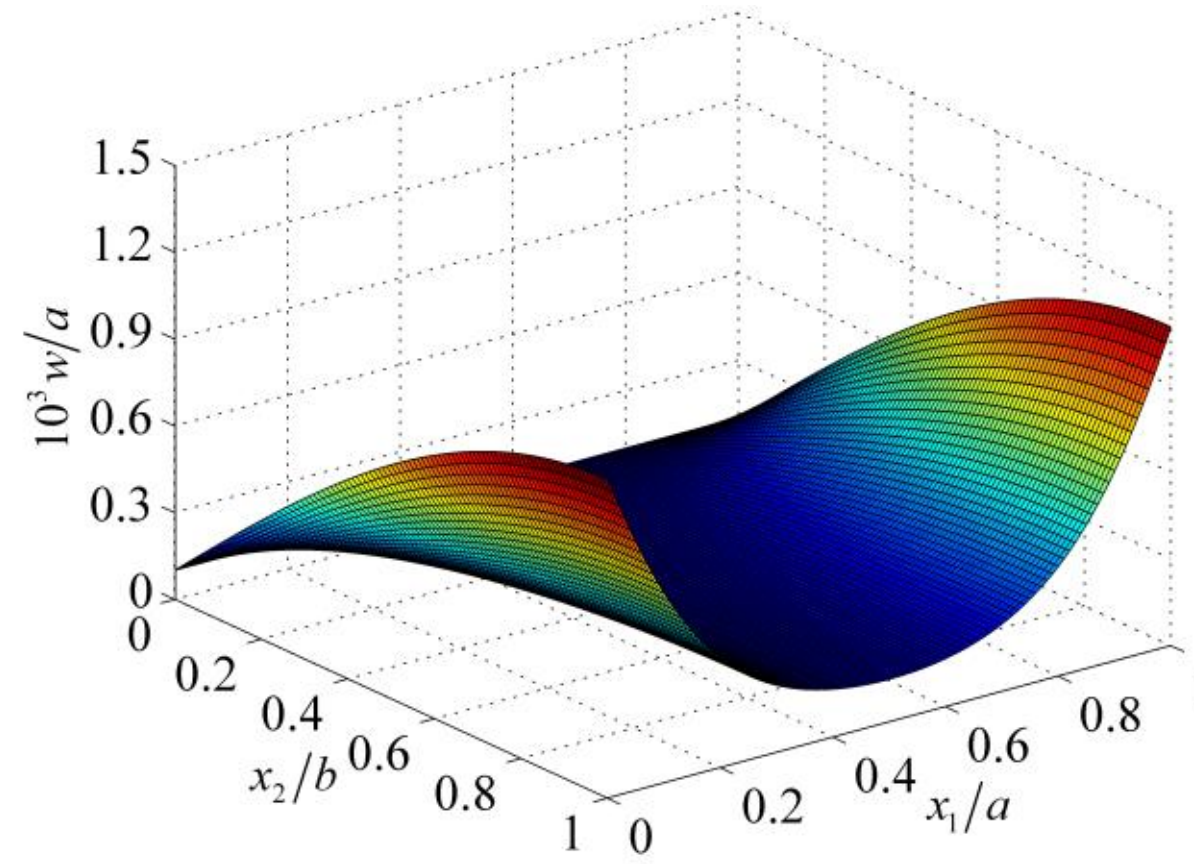

(a)

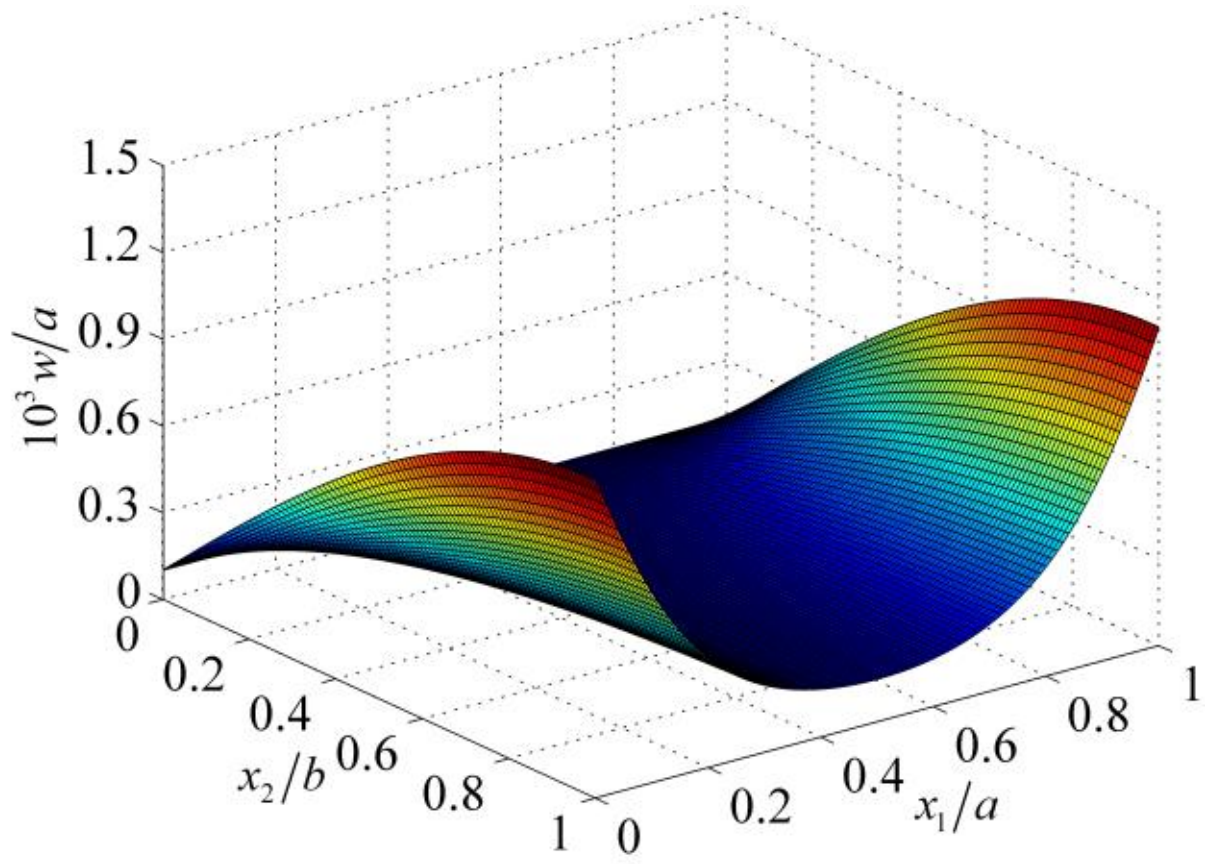

(b)

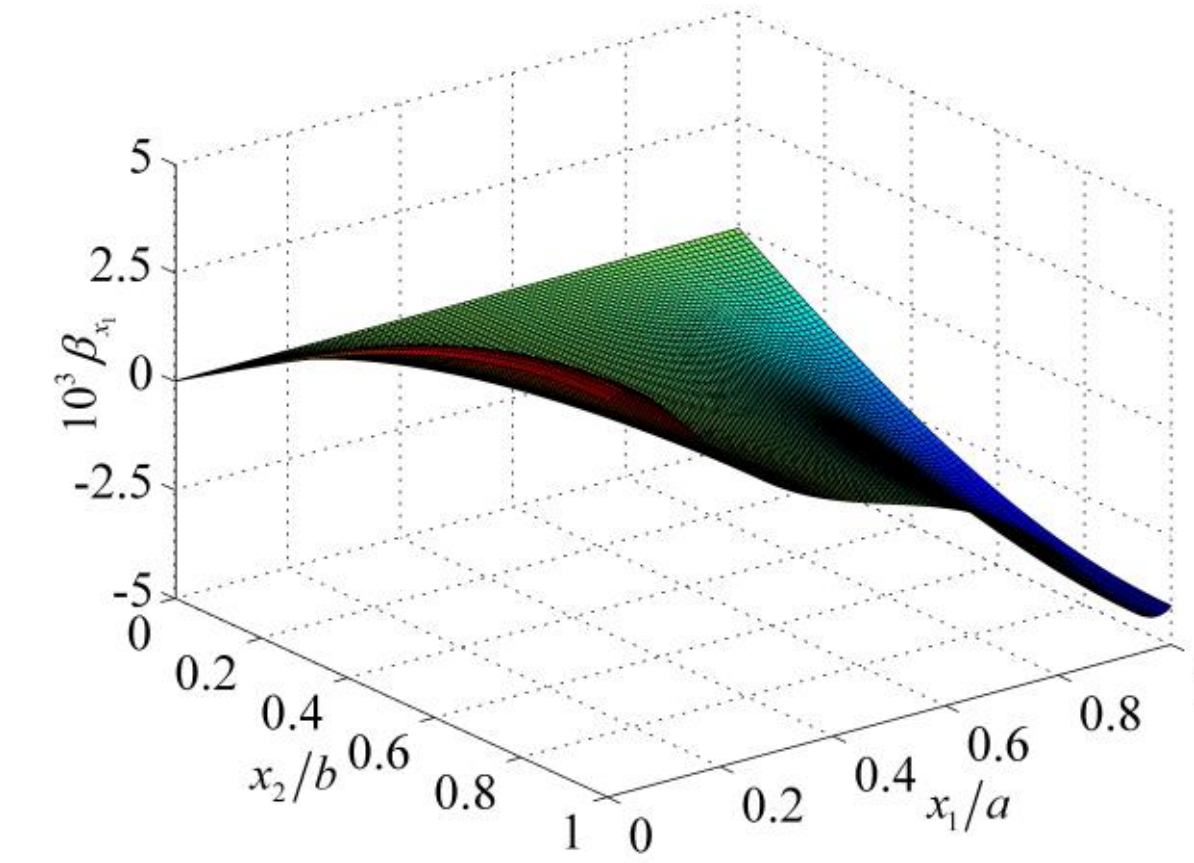

(c)

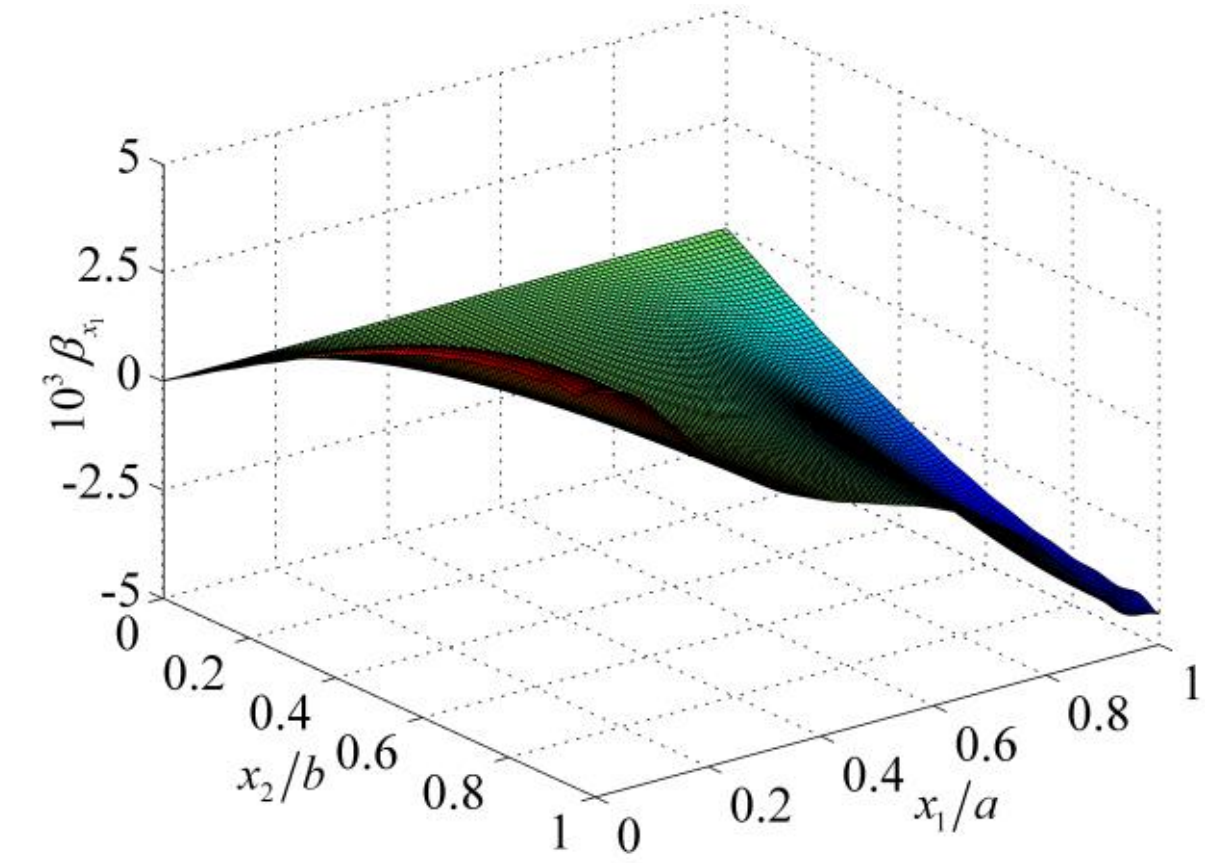

(d)

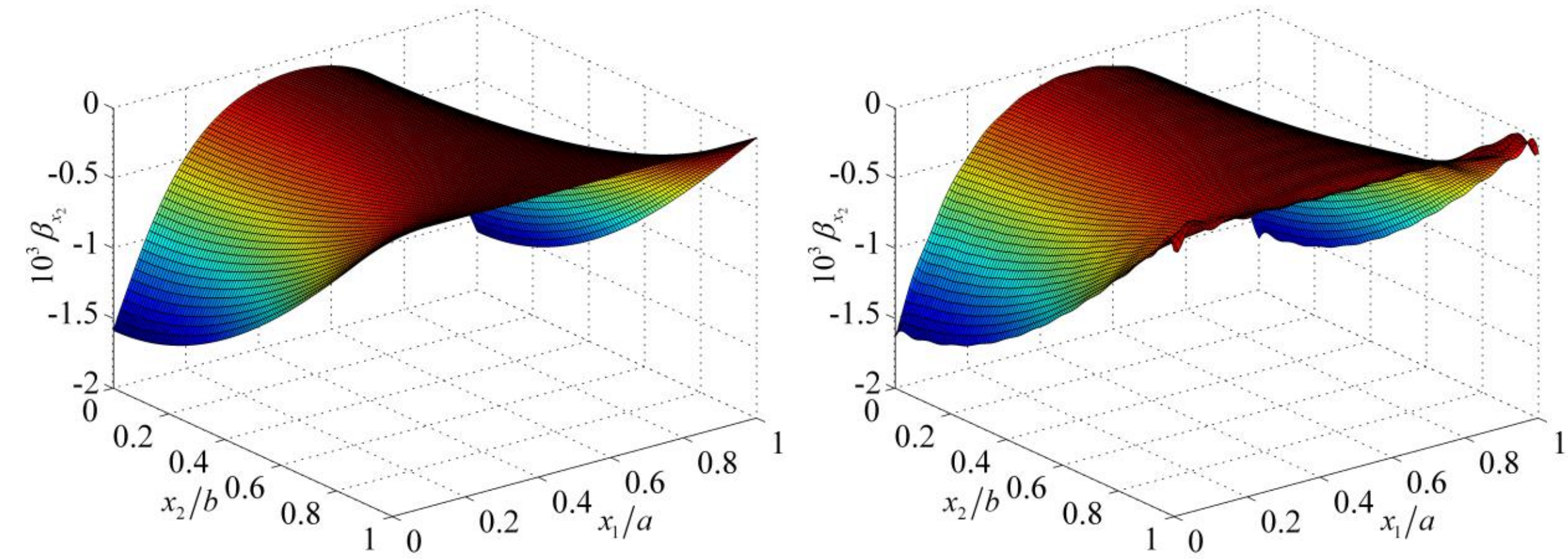

(e) (f)

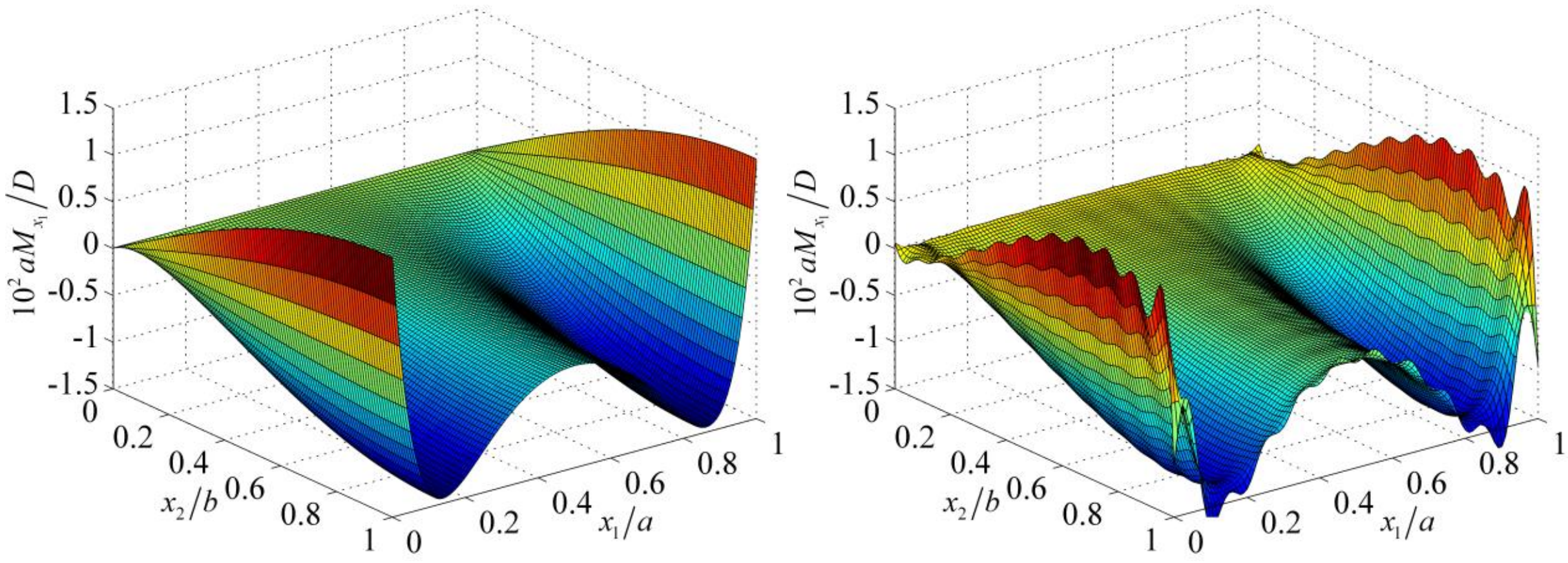

(g) (h)

Figure 9: The deflection, rotation and moment distributions for $k_r = 10^4$ and $G_{pr} = 190$:
(a) $w$ (reference), (b) $w$ (calculated), (c) $\beta_{x_1}$ (reference), (d) $\beta_{x_1}$ (calculated),
(e) $\beta_{x_2}$ (reference), (f) $\beta_{x_2}$ (calculated), (g) $M_{x_1}$ (reference), (h) $M_{x_1}$ (calculated).

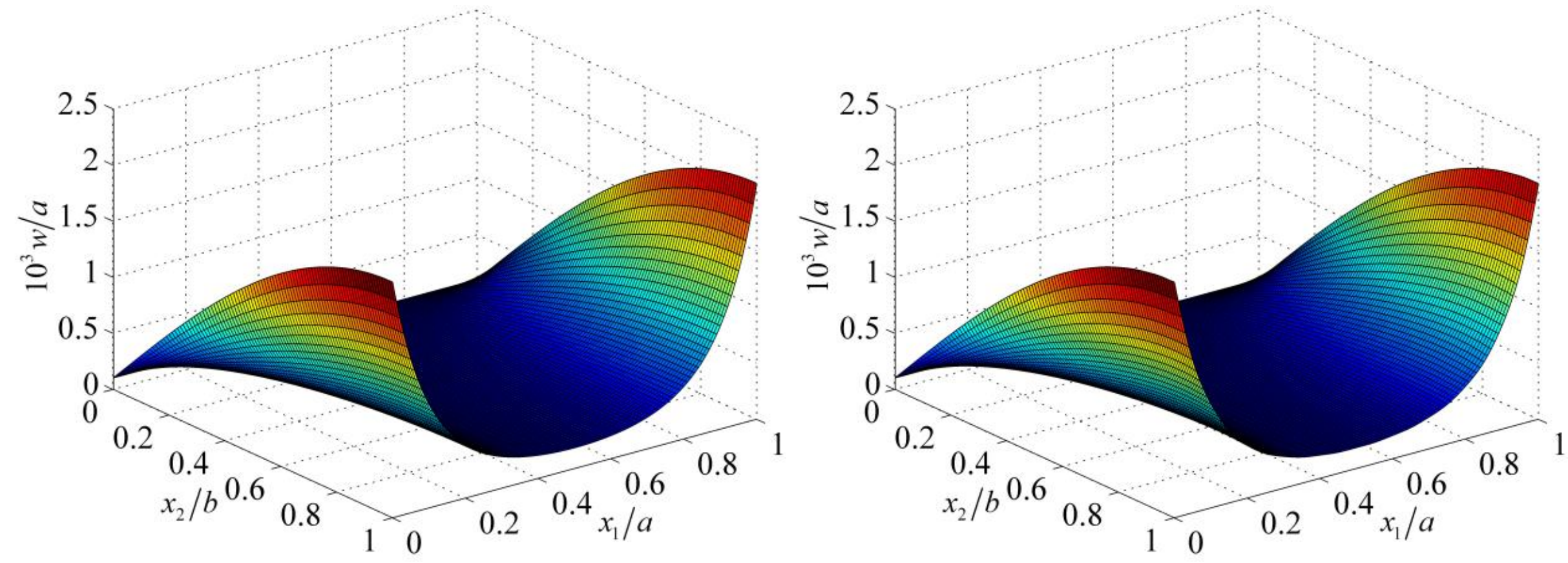

(a) (b)

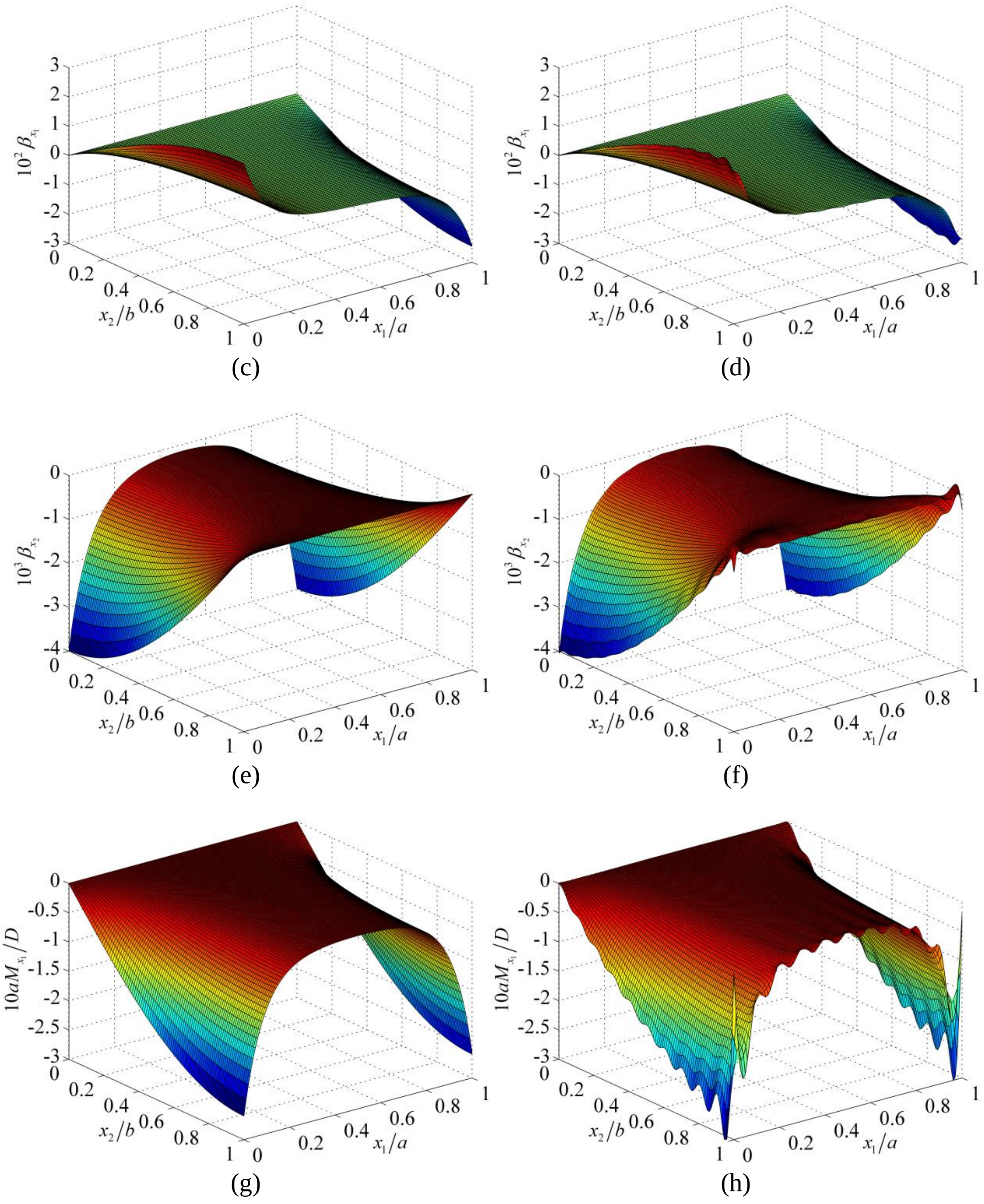

Figure 10: The deflection, rotation and moment distributions for $k_r = 10^4$ and $G_{pr} = 300$:
(a) $w$ (reference), (b) $w$ (calculated), (c) $\beta_{x_1}$ (reference), (d) $\beta_{x_1}$ (calculated),
(e) $\beta_{x_2}$ (reference), (f) $\beta_{x_2}$ (calculated), (g) $M_{x_1}$ (reference), (h) $M_{x_1}$ (calculated).

## 6. Conclusions

Plate on elastic foundations is a common model for several types of engineering structures. In this paper, the usual structural analysis of a Reissner plate on the Pasternak foundation is extended to a multiscale analysis of a system of linear differential equations with general boundary conditions and a wide spectrum of model parameters. It is concluded that:

1. We derive the Fourier series multiscale solution of the bending problem of the Reissner plate on the Pasternak foundation.
2. We employ the minimum potential energy method, rather than the usually used Fourier coefficient comparison method, to derive the final system of equations to be solved for the Fourier series multiscale solution.
3. We investigate the convergence characteristics of the obtained Fourier series multiscale solution.
4. We reveal the multiscale characteristics of the bending problem of the Reissner plate on the Pasternak foundation.

The preliminary study on application to bending of the Reissner plate on the Pasternak foundation verifies the effectiveness of the present Fourier series multiscale method, and the Fourier series multiscale solution potentially offers a powerful means for the analysis of some new types of concrete (such as FRP-reinforced concrete or ultra-high-performance concrete) pavements of engineering structures.